\input amstex
\magnification 1200

\loadmsbm
\parindent 0 cm

\define\nl{\bigskip\item{}}
\define\snl{\smallskip\item{}}
\define\inspr #1{\parindent=20pt\bigskip\bf\item{#1}}
\define\iinspr #1{\parindent=27pt\bigskip\bf\item{#1}}
\define\einspr{\parindent=0cm\bigskip}

\define\ot{\otimes}
\define\R{\Bbb R}
\define\C{\Bbb C}
\define\Z{\Bbb Z}
\define\T{\Bbb T}
\define\N{\Bbb N}

\input amssym


\centerline{\bf The Haar measure on some locally compact}
\centerline{\bf quantum groups}
\bigskip
\centerline{\it Alfons Van Daele} 
\bigskip\bigskip\bigskip
{\bf Abstract}
\nl
A locally compact quantum group is a pair $(A,\Phi)$ of a C$^*$-algebra 
$A$ and a $^*$-homomorphism $\Phi$ from $A$ to the multiplier algebra 
$M(A\ot A)$ of the minimal C$^*$-tensor product $A\ot A$ satisfying 
certain assumptions (see [K-V1] and [K-V2]). One of the assumptions is the 
existence of the Haar weights. These are densely defined, lower 
semi-continuous faithful KMS-weights satisfying the correct invariance 
properties.
\snl
Many examples of C$^*$-algebras with a comultiplication arise from 
quantizations of classical locally compact groups. These are first 
obtained on the Hopf $^*$-algebra level and then lifted to the 
C$^*$-algebra context. This step is usually rather complicated (but 
interesting analysis is involved). It is necessary if one wants to have 
the Haar weights. It is Woronowicz who has done remarkable work in this 
direction. However, his technique to pass from the Hopf $^*$-algebra 
level to the C$^*$-level does not give the Haar weights. In this paper, 
we will study a recent example of Woronowicz, the quantum $az+b$-group, 
and obtain the Haar weights. We use a technique that is  
useful in other cases as well (as we in fact show at the 
end of our paper). 
\snl
An important feature of the example we study here is that the Haar 
weights are not invariant, but only relatively invariant with respect 
to the scaling group (coming from the polar decomposition of the 
antipode). This phenomenon was expected from the theory, but up to now, 
no such example existed. It is another indication that the notion of a 
locally compact quantum group, as introduced and studied by Kustermans 
and Vaes in [K-V1] and [K-V2], is the correct one.    
\nl
March 2001 ({\it Preliminary version})
\nl
\vskip 2.3 cm
\hrule
\medskip
Department of Mathematics, K.U.\ Leuven, Celestijnenlaan 200B,
B-3001 Heverlee (Belgium). 

E-mail address : Alfons.VanDaele\@wis.kuleuven.ac.be

\newpage


\bf 1. Introduction \rm
\nl
Let $G$ be a locally compact group. Denote by $C_0(G)$ the C$^*$-algebra  
of continuous complex functions on $G$, tending to $0$ at infinity. 
Identify $C_0(G\times G)$ with the C$^*$-tensor product $C_0(G)\ot 
C_0(G)$. Also identify the C$^*$-algebra $C_b(G\times G)$ of bounded 
continuous complex functions on $G\times G$ with the multiplier algebra 
$M(C_0(G)\ot C_0(G))$ of $C_0(G)\otimes C_0(G)$. Then the product in $G$ 
gives rise to a non-degenerate $^*$-homomorphism $\Phi:C_0(G) \to  
M(C_0(G)\ot C_0(G))$ given by $\Phi(f)(p,q)=f(pq)$ whenever $p,q\in G$ 
and $f\in C_0(G)$. This $\Phi$ is called a comultiplication. The 
inverse in the group gives the antipode $S$ mapping $C_0(G)$ to itself and 
is defined by $(Sf)(p)=f(p^{-1})$. Finally, the left and right Haar 
measures induce weights $\varphi$ and $\psi$ on $C_0(G)$, defined by 
integration with respect to these measures. Left and right invariance 
can be expressed (using $\iota$ for the identity map) 
as $(\varphi\ot\iota)\Phi(f)=\varphi(f)1$, respectively 
$(\iota\ot\psi)\Phi(f)=\psi(f)1$ for $f\geq 0$ and the correct 
interpretation of these formulas.
\snl
The above passage from the locally compact group to the C$^*$-algebra 
with comultiplication, antipode and invariant weights, provides the 
basic idea for the development of the theory of locally compact {\it 
quantum} groups. Roughly speaking, a locally compact quantum group is a 
pair $(A, \Phi)$ of a C$^*$-algebra $A$ and a $^*$-homomorphism  
$\Phi:A\to M(A\ot A)$ of $A$ to the multiplier algebra of the minimal 
C$^*$-tensor product $A\ot A$ satisfying certain (natural) properties 
and admitting nice left and right invariant weights. We refer to [K-V2] 
for the precise notion and the development of the concept.
\snl
In this paper, we deal with a class of {\it examples} of such locally compact 
quantum groups.
\snl
For understanding the importance of these examples in the development 
of the theory, we have to mention some historical facts about some of 
the intermediate steps preceding the general theory of locally compact 
quantum groups as it is known now.
\snl
The origin lies in an attempt to generalize the Pontryagin duality 
theorem for locally compact abelian groups. There were many 
intermediate steps that eventually led to the theory of Kac algebras 
where the duality was again obtained within the same category (as is 
the case for abelian groups). We refer to [E-S] for the theory of Kac 
algebras. Also in the introduction of [E-S], more can be found on the 
history between the original theorem of Pontryagin and the Kac 
algebras.
\snl
Apart from being  rather complicated, the Kac algebra theory 
turned out to be unsatisfactory from different points of view. The 
framework is von Neumann algebra theory while C$^*$-algebras are more 
natural since we are quantizing a topological group. The antipode $S$ 
was assumed to be a $^*$-anti-homomorphism (or equivalently satisfying
$S^2=\iota$). It was known for some time that there were Hopf algebras 
not having this property. Hence, it should be no surprise that this 
restriction on the antipode was violated in the examples of Drinfel'd 
and Jimbo (in quantum group theory) ([D] and [J]).
\snl
Probably, E. Kirchberg was the first to propose a generalization of the 
Kac algebras without this restriction on the antipode [Ki]. He proposed some 
kind of {\it polar decomposition} of the antipode involving a so-called 
unitary antipode and a one-parameter group of $^*$-automorphisms (the 
scaling group). Such a theory was developed by Masuda and Nakagami 
[M-N]. Still, their axioms are rather complicated and the theory is 
formulated in the von Neumann algebra setting. Now, Masuda, Nakagami 
and Woronowicz are reformulating this theory in the C$^*$-algebra 
framework. Their work has been announced on several occasions (e.g.\ 
[M-N-W]) but it is not yet available.
\snl
At the same time, we developed the theory of multiplier Hopf 
$^*$-algebras with positive integrals (see [VD2] and [VD3]). This led to 
a certain class of locally compact quantum groups ([K1] and [K-VD]). 
One of the 
remarkable features was that the Haar measures where not necessarily 
invariant, but 
only relatively invariant with respect to the scaling group. On the 
other hand, in the theory of Masuda, Nakagami and Woronowicz, the Haar 
measures are assumed to be invariant for the scaling group. This of 
course rose the question whether or not examples existed where this was 
not the case. The more recent development of the theory of locally 
compact quantum groups by Kustermans and Vaes ([K-V2]) provided another 
strong argument for the possibility of non-invariance. Indeed, they 
start from a rather simple and natural set of axioms and arrive at Haar 
weights, only relatively invariant with respect to the scaling group.
\snl
While this development of the theory was going on, Woronowicz provided 
us with many examples, in general rather complicated ones, but 
involving interesting analysis. He developed these examples from a 
certain point of view with typical techniques but these do not give the 
Haar measure. An example of this is the quantum $E(2)$ and its dual. 
The Haar measure here was obtained by S. Baaj in [B1] and [B2]. Of 
course, the Haar measure is very important as it is, in a way, the main 
reason for studying quantum groups in the operator algebra framework. 
\snl
In a separate paper, we plan to give an approach to some of these 
examples, starting from dual pairs of Hopf $^*$-algebras (see [VD4]). 
In the present paper, we will 
fully use the results as obtained by Woronowicz. And although the 
$ax+b$-group was the first we considered, we will only treat this case  
at the end of the paper. The reason is that Woronowicz, more recently, 
developed another example that turned out to have the same phenomenon 
of non-invariance and is more easy (one of the reasons being that this 
example also exists on the pure Hopf $^*$-algebra level). It is a 
quantization of the $az+b$-group.
\snl
So, in this paper, we will describe the quantum $az+b$-group as 
introduced by Woronowicz in [W5], we will construct the Haar measures 
and we will see that they are not invariant but only relatively 
invariant with respect to the scaling group.
\snl
This paper is not only important because it provides the first example 
having this property of non-invariance. The techniques that we use to 
obtain the Haar measures and to prove the invariance are also new. And 
they are expected to work in many more cases than the examples we 
describe here. In fact, we strongly believe that the techniques used here 
can contribute to a possible general existence theorem for the Haar 
measure on a general locally compact quantum group. Recall that the 
existing theories all assume the existence of the Haar measures (except 
for special cases - like the compact and discrete quantum groups). In a 
forthcoming pair of papers [VD5] and [VD6], 
we plan to develop these ideas further. In 
the first one, we will solely look at the algebraic aspects while in 
second one, we will 
treat the C$^*$-algebra and von Neumann algebra versions.
\snl
Before giving the content of this paper, let us say something about the 
style of it. In general, the theory of locally compact quantum groups 
is not very easy. The examples are even more complicated (although some 
of the more recent ones, as obtained by Vaes and Vaynerman (see [V-V] 
and [V]) seem to be simpler). Nevertheless, the theory is rich and very 
beautiful, and is expected to have nice applications. In this paper, we 
treat examples, but this is done in such a way that also a better 
understanding of the general theory can be obtained.
\snl
Moreover, this work here is based very heavily on the work of 
Woronowicz. To make the paper still, up to some degree, self-contained, 
we will recall that part of the paper of Woronowicz that will be needed 
to understand this paper. In fact, we will take a slightly different 
point of view with respect to manageability and the relation with the 
polar decomposition of the antipode. This is another reason why we 
want to include some of the results of [W5] so as to formulate them in 
a way suitable for us. We will also, where more convenient, use the von 
Neumann algebra framework rather than the C$^*$-algebra setting. This 
is sometimes easier (although further away from the original 
formulation). See e.g.\ section 4 where we use von Neumann algebra 
theory to construct the Haar measure.
Finally, as we mentioned already, the ideas we 
are using here are expected to have more general applications and so, 
we will explain these in a fairly great detail.
\snl
On the technical level, we have chosen to be as precise as possible. 
Sometimes though, we have prefered a more loose way of writing in order 
not to let the technical details prevent a better understanding of what 
is really going on. In these cases, the reader should always be able to 
make the arguments completely rigourous.
\nl 
Now we come to the {\it structure} and the {\it content} of this paper.
\snl
At the end of this introduction, we will collect some basic notions and 
give some standard references. In section 2 of the paper, we give the 
quantum $az+b$-group on the Hopf $^*$-algebra level. Although this part 
is not really necessary for the rest of the paper, it is rather easy 
and very instructive. It will help to understand the formulas in the 
further sections. The polar decomposition of the antipode is one of the 
results that can already be obtained in this purely algebraic context.
\snl
In section 3, we describe the C$^*$-algebra  and the 
comultiplication. We start from the formulas in section 2 and realize 
the generators of the Hopf $^*$-algebra with operators on a Hilbert 
space, having the correct (strong versions of the) commutation rules. 
The C$^*$-algebra is obtained by taking suitable functions of the 
generators. The comultiplication is given using the implementation by a 
suitable multiplicative unitary. For this multiplicative unitary, we 
need to refer to [W5]. The existence of the comultiplication is highly 
non-trivial and the work of Woronowicz is very important here. Again, 
at the end of the section, we treat the antipode and its polar 
decomposition. All of this is due to Woronowicz. But, as we mentioned 
before, we take a slightly different point of view with respect to 
manageability - a notion we try to avoid - and the relation with the polar 
decomposition of the antipode. In any case, we will indicate clearly 
what these small differences are and explain the relation of the 
different multiplicative unitaries that we encounter. We also say more 
about our point of view in section 6 where we give some conclusions and 
discuss some perspectives. 
\snl
In section 4, we start with the construction of the weight that is the 
candidate for the right Haar measure (see theorem 4.4). 
Before we prove that it is really 
right invariant, we explain why this result can be expected. The idea 
behind this argument will be seen, in section 5, to be crucial, also 
for other examples. And in section 6, we will explain why this idea 
might be a basis for a possible general existence theorem for the Haar 
measures on locally compact quantum groups. After giving this 
information, we prove the invariance and we show that we get a locally 
compact quantum group in the sense of Kustermans and Vaes. We prove 
that the Haar weights are not invariant but only relatively invariant 
with respect to the scaling group. Finally we give an explicit formula 
for the right regular representation. It gives another multiplicative unitary 
and we discuss the relation of this new one with the original one as 
found in Woronowicz' work. 
\snl
In section 5 we briefly treat other examples of this type : 
the quantum $az+b$-group (with real deformation parameter) and the quantum 
$ax+b$-group. We will also mention the quantum $E(2)$ and its 
dual. We will not go into details here as the method is 
completely the same as for the quantum $az+b$-group described in the 
previous sections. We will give the formulas for the Haar measure and 
indicate the main points (and possible, more fundamental differences
with the $az+b$-case).
\snl
Finally, in section 6, we will draw some conclusions and discuss 
possible further research, based on the ideas and techniques introduced 
in this paper.
\snl
At the end of the paper, in an appendix, we recall some of the 
well-known facts about different forms of the Heisenberg commutation 
relations. These relations, in various cases, appear here and there in 
the examples that we study in this paper.
\nl
Now, let us recall some {\it basic notions} and give {\it standard 
references}.
\snl
We will work with $^*$-algebras over the complex numbers, with or 
without identity $1$. An element $a$ in a $^*$-algebra is called normal 
if $aa^*=a^*a$ and unitary if moreover $aa^*=a^*a=1$. It is called 
self-adjoint if $a=a^*$. When we are dealing with bounded operators on 
a Hilbert space $\Cal H$, the same terminology is used by considering 
the $^*$-algebra $\Cal B(\Cal H)$ of all bounded operators on $\Cal H$ 
with the adjoint as the involution. Of course, for unbounded 
operators, the notion of normal and self-adjoint operators is more 
special. We refer to [K-R] for the theory of (unbounded) operators on a 
Hilbert space.
\snl
In section 2, where we consider the Hopf algebra level, we will be 
dealing with algebraic tensor products. In the rest of the paper, we 
will have topological (i.e.\ completed) tensor products (of Hilbert 
spaces, C$^*$-algebras, von Neumann algebras, ...). We will in all 
cases use the symbol $\ot$ but in general, 
it should be clear from the context what 
we have.
\snl
Recall the definition of a Hopf algebra. It is a pair $(H,\Delta)$ of 
an algebra $H$ over $\Bbb C$ with an identity $1$ and a unital 
homomorphism $\Delta : H \to H\ot H$ satisfying coassociativity 
$(\Delta\ot\iota)\Delta=(\iota\ot\Delta)\Delta$ (recall that we use 
$\iota$ to denote the identity map) and such that there is 
a counit $\varepsilon$ and an antipode $S$. The counit is characterized 
as a linear map $\varepsilon : H \to \Bbb C$ such that $(\varepsilon\ot 
\iota)\Delta(a)=(\iota\ot\varepsilon)\Delta(a)=a$ for all $a\in H$. It 
is uniquely determined by this property and it is a homomorphism. The 
antipode is a linear map $S:H\to H$, characterized by the equations
$m(S\ot\iota)\Delta(a)=m(\iota\ot S)\Delta(a)=\varepsilon(a)1$ for all 
$a$ in $H$ (where $m$ denotes the multiplication as a map from $H\ot H$ 
to $H$). Again $S$ is unique but it is a anti-homomorphism. When $H$ is 
a $^*$-algebra, then $(H,\Delta)$ is called a Hopf $^*$-algebra if 
moreover $\Delta$ is a $^*$-homomorphism. In this case, $\varepsilon$ 
is also a $^*$-homomorphism while $S(S(a)^*)^*=a$ for all $a\in H$. In 
fact, $S$ is a $^*$-anti-homomorphism if and only if $S^2=\iota$. We 
refer to [A] and [S] for the theory of Hopf algebras and to [VD1] for the 
theory of Hopf $^*$-algebras and dual pairs of Hopf $^*$-algebras.
\snl
We will use $\Delta$ for the comultiplication in Hopf algebras and 
$\Phi$ for a comultiplication on a C$^*$-algebra $A$ (where $\Phi$ is a 
$^*$-homomorphism  from $A$ to the multiplier algebra $M(A\ot A)$ of 
the minimal C$^*$-tensor product $A\ot A$). We will also use $\Phi$ for 
a comultiplication on a von Neumann algebra $M$ (where now $\Phi$ maps 
into the von Neumann tensor product $M\ot M$). We refer to [K-R], [S-Z], 
[Sa] and [P] for the 
basics of C$^*$-algebra  and von Neumann algebra theory.
\snl
We will use the theory of left Hilbert algebras to construct the Haar 
weight. We refer to [St] for the theory of Hilbert algebras and 
for the theory of weights on C$^*$-algebras and von Neumann algebras. 
The standard procedure to construct a weight from a left Hilbert 
algebra can also be found in [St]. We will also use the standard notations 
related to weights. When $\psi$ is a weight on a C$^*$-algebra $A$, 
then we use $\Cal N_\psi$ for the left ideal of elements $x$ in $A$ 
satisfying $\psi(x^*x)<\infty$ and $\Cal M_\psi$ for the subalgebra 
spanned by the positive elements $x$ in $A$ satisfying 
$\psi(x)<\infty$. 
\nl
\nl 
\bf Acknowledgements \rm
\snl
First of all, I would like to thank S.L.\ Woronowicz for providing me 
with the preliminary versions of his work on the quantum $ax+b$- and 
$az+b$-groups. I have profited very much from many discussions on this 
topic with him on various occasions.
\snl
Most of the work, that ended up in this article, was done while 
visiting the University of Trondheim and I would like to thank my 
colleagues there, in particular, M.\ Landstad and C.\ Skau for the 
hospitality.
\snl
I also thank K.\ Schm\"udgen for giving me the opportunity to talk 
about this work at a meeting in Bayrischzell.
\nl
\nl


\bf 2. The Hopf $^*$-algebra \rm
\nl
In this section, we will define the Hopf $^*$-algebra that is the basis 
for 
the example of the $az+b$-group. 
We will introduce it via some intermediate steps. These 
various Hopf algebras play also a role in the other examples that we 
plan to treat briefly in section 5. The reason why we precisely picked 
this example to treat in detail is because it exists already at the 
Hopf $^*$-algebra level. However, it should be mentioned that the 
results about the Hopf algebra are not needed to treat the 
C$^*$-algebra case. But of course, this is highly instructive and, 
following the spirit of the paper, we include these results here. They 
will yield a better understanding of the main part of the paper.
\snl
We start with the following well-known Hopf algebra. It is probably the 
simplest example of a non-commutative, non-cocommutative Hopf algebra. 
It is a deformation of the Hopf algebra of polynomial functions on the 
classical $ax+b$-group.

\inspr{2.1} Proposition \rm
Let $\lambda$ be any non-zero complex number.  Let $H$ be the
algebra over $\Bbb C$ with identity generated by two elements $a$ and $b$ such
that $a$ is invertible and $ab = \lambda ba$.  There is a comultiplication
$\Delta$ on $H$ defined by
$$\align
\Delta(a) & = a \otimes a \\
\Delta(b) & = a \otimes b + b \otimes 1.
\endalign
$$
The pair $(H, \Delta)$ is a Hopf algebra. The counit $\varepsilon$ and antipode
$S$ are given by
$$\alignat 2
\varepsilon(a) &= 1 & \qquad \qquad \quad S(a) &= a^{-1} \\
\varepsilon(b) &= 0 & \qquad \qquad \quad S(b) &= -a^{-1}b.
\endalignat$$
\einspr

The proof of this result is very simple and straightforward. A possible 
reference is [VD1] where similar examples are worked out in more detail. 
The methods can also be applied here.
\snl
This Hopf algebra can be made into a Hopf $^*$-algebra by imposing the 
conditions that $a$ and $b$ are self-adjoint elements. In this case, 
one should require that $|\lambda|=1$ for the commutation relation to 
be compatible with the $^*$-operation. Also this result is well-known, 
again, for an example of this type, see e.g.\  example 2.6 in [VD1]. 
\snl
This is a simple, non-trivial example of a Hopf $^*$-algebra. In this 
paper, we are interested in the C$^*$-versions of these Hopf 
$^*$-algebra structures. This example, though very simple on the Hopf 
$^*$-algebra level, becomes a lot more complicated on the C$^*$-level. 
The problem is not so much caused by the $^*$-algebra structure as this 
can easily be lifted to the C$^*$-framework. The real problem is 
caused by the comultiplication (see e.g.\ [W-Z]). We will 
discuss this example briefly in section 5.
\snl
There is however another way to produce a Hopf $^*$-algebra from this 
Hopf algebra. For this example, the passage to the C$^*$-level is possible, 
although still quite involved (as we will see in the next section), 
 without adding some extra generator as one has to do to obtain 
the $ax+b$-group on the quantum level in the framework of locally 
compact quantum groups.
\snl 
In stead of putting a $^*$-algebra structure on the Hopf algebra of 
proposition 2.1, we now take a direct sum of two copies of this Hopf 
algebra and define the involution by sending the generators of one copy 
to the corresponding generators of the other copy. This results in the 
following.
\inspr{2.2} Proposition \rm  
Let $\lambda$ be a non-zero complex number.  Let $H$ be the
$^*$-algebra over $\Bbb C$ with identity generated by normal 
elements $a$ and $b$
such that $a$ is invertible,  $ab = \lambda ba$ and $ab^* = b^*a$.  
There is a comultiplication $\Delta$ defined on $H$ as in proposition 
2.1. The pair $(H,\Delta)$ is a Hopf $^*$-algebra.
The counit and the antipode are given by the same formulas as in 2.1. 
Now also 
$$\align S(a^*)&=(a^*)^{-1} \\ S(b^*)&=-(a^*)^{-1}b^*. \endalign$$
\einspr

There is no need for a restriction
of the number $\lambda$.  By taking the adjoint of $ab = \lambda ba$ we get 
$a^* b^* = \mu b^* a^*$ where $\mu = \overline{\lambda^{-1}}$.  
And by taking adjoints in
$ab^* = b^* a$, we get $a^* b = ba^*$.  And we indeed have the
direct sum of two copies of the Hopf algebra in proposition 2.1. (possibly with
different $\lambda$-factor).
\snl
It is well known that the Hopf algebra in 2.1 can be paired with itself.  
Such a pairing is given by 
$$\alignat 2
\langle a,a \rangle &= \lambda &\qquad\qquad \langle b,a \rangle &= 0 \\
\langle a,b \rangle &= 0   &\qquad\qquad \langle b,b \rangle &= t
\endalignat$$
where $t$ is any complex number. Again see [VD1] for a similar case.
\snl
If we endow this Hopf algebra with the $^*$-structure making $a$ and 
$b$ self-adjoint (and assuming $|\lambda|=1$), we get a dual pair of 
Hopf $^*$-algebras (i.e.\ such that also $\langle x^*,y\rangle=\langle 
x,S(y)^*\rangle^-$ for all $x,y\in H$), if $t$ is purely imaginary.
\snl
Considering the direct sum of two copies of this Hopf algebra with the 
$^*$-algebra structure as in 2.2, we still get a self-pairing between 
Hopf $^*$-algebras. The pairing respects the direct sum structure in 
the sense that the algebra generated by $a$ and $b$ is paired with 
itself as above but has trivial pairing with the algebra generated by 
$a^*$ and $b^*$. See [VD4] for more details on this approach to these 
examples. This pairing gives a better understanding of the formula for 
the multipicative unitary that we will have in the next section 
(definition 3.7).
\snl
The pairing that we describe above will be degenerate when $\lambda$ is 
a root of unity. In that case some power of $a$ will be $1$ and some 
power of $b$ will be $0$ (in the quotient).
\snl
In the case of the Hopf $^*$-algebra given in proposition 2.2, there is 
another quotient. This quotient turns out to be still self-dual (again 
see [VD4]). And it is precisely this Hopf $^*$-algebra that will be the 
one studied here. We obtain it in the following proposition.

\inspr{2.3} Proposition \rm
Let $n$ be a non-zero natural number and put $\lambda=\exp\frac{2\pi 
i}{n}$. Let $H$ be the $^*$-algebra over $\Bbb C$ with identity 
generated by two normal elements $a$ and $b$ such that $a$ is 
invertible, $ab=\lambda ba$, $a^*b=ba^*$ and moreover such that $a^n$ 
and  $b^n$ are self-adjoint. There is a comultiplication $\Delta$ on 
$H$ making $H$ into a Hopf $^*$-algebra, given by the formulas in 2.1. 
The counit and antipode are also given by the same formulas (as in 2.1 
and 2.2).
\snl
\bf Proof : \rm To prove that $\Delta$, defined as in 2.1, is 
well-defined, we must show that $(a \otimes a)^n$ and 
$(a \otimes b + b \otimes 1)^n$ are self-adjoint in $H\ot H$. This is 
obvious in the first case. In the second case, as $ab=\lambda ba$, we 
can write 
$$(a \otimes b + b \otimes 1)^n=\sum_{k=0}^{n}c_k^n b^{n-k}a^{k}\ot b^k 
$$
where the coefficients $c_k^n$ are complex numbers, depending on 
$\lambda$.
Using standard techniques (see e.g.\ [K-S]), 
it follows from the fact that $n$ is the 
smallest number such that $\lambda^n=1$, that
$$(a \otimes b + b \otimes 1)^n=a^n \otimes b^n + b^n \otimes 1$$    
and this is indeed self-adjoint.
\snl
To prove that we still have a Hopf $^*$-algebra, it is also necessary 
to show that the extra conditions are compatible with the counit and 
the antipode. The counit gives no problem and also the behaviour of $S$ 
on $a$ represents no difficulty.  It remains to verify that $S(b^n)^*= 
S^{-1}(b^n)$.
\snl
Now, a straightforward calculation, using standard techniques, shows 
that
$$S(b^n)=(-1)^n\lambda^{\frac12 n(n-1)}a^{-n}b^n.$$
On the other hand, $S^{-1}(b)=-ba^{-1}$ so that
$$S^{-1}(b^n)=(-1)^n\lambda^{-\frac12 n(n-1)}b^na^{-n}.$$
The  equation to verify now follows from the fact that $|\lambda|=1$.
\einspr

In this proposition, one should really take $n\geq 2$. However, the 
result still is true for $n=1$  but then, it is trivial. In that case, 
the algebra is abelian and the generators $a$ and $b$ are self-adjoint.
\snl
Later, the extra algebraic conditions on $a$ and $b$ will be formulated 
by a spectral condition on the operators that we will take for $a$ and 
$b$ (see definition 3.1). Thus, in the C$^*$-algebra context, the quotient in 2.4 of 
the example in 2.3 is obtained by imposing spectral conditions. This is 
a known phenomenon (see the work on the quantum $E(2)$ ([W4]), on the 
quantum $ax+b$-group ([W-Z]), ...). Also the need to pass to a discrete 
deformation parameter is no longer considered to be peculiar. It is 
just part of the quantization process.
\snl
Now, we will obtain some more information about the antipode. In fact, 
we have a {\it polar decompostion} for the antipode. That such a polar 
decomposition already exists on the Hopf $^*$-algebra level, is not so 
uncommon. See e.g.\ [K2] where it is proved that this is always the 
case for multiplier Hopf $^*$-algebras with positive integrals, in 
particular for discrete and compact quantum groups.

\inspr{2.4} Proposition \rm 
Let $(H,\Delta)$ be the Hopf $^*$-algebra obtained in 
proposition 2.3. There is an involutive $^*$-anti-automorphism $R$ of 
$H$ and a one-parameter group $\{\tau_t\mid t\in\Bbb R\}$ of 
$^*$-automorphisms such that $t\to \tau_t(x)$ is analytic for all $x\in 
H$ and such that $S=R\tau_{-\frac{i}{2}}=\tau_{-\frac{i}{2}}R$ (where 
$\tau_z$ is the analytic extension of $\tau_t$ to $z\in \Bbb C$).
For all $t$, we have that $R$ and $\tau_t$ commute. Moreover
$$\Delta(R(x))=\sigma(R\ot R)\Delta(x)$$
where $\sigma$ denotes the flip from $H\ot H$ to itself given by 
$\sigma(x\ot y)=y\ot x.$ Also
$$\Delta(\tau_t(x))=(\tau_t\ot \tau_t)\Delta(x)$$
for all $x\in H$.
\snl
\bf Proof: \rm
As it is easier to define $\tau_t$ than to define $R$, we start with this 
one-parameter group. Let $\tau_t(a)=a$ and $\tau_t(b)=e^{\frac{2\pi 
t}{n}}b$ for all $t\in \Bbb R$. It is easy to verify that these 
formulas yield a one-parameter group of $^*$-automorphisms of the Hopf 
$^*$-algebra.
\snl
This is an analytic one-parameter group in the sense that the map $t\to 
f(\tau_t(x))$ is analytic for all $x\in H$ and all linear functionals 
$f$ on $H$.
\snl
It also satisfies $\Delta(\tau_t(x))=(\tau_t\ot \tau_t)\Delta(x)$ for 
all $x\in H$ and $t\in \Bbb R$. Finally observe that $\tau_t$ commutes 
with $S$.
\snl
Next comes the definition of $R$. By analyticity, we can define the 
automorphism $\tau_{\frac{i}{2}}$ and we let $R=S\tau_{\frac{i}{2}}$. 
Because $S$ commutes with $\tau_t$, it will also commute with 
$\tau_{\frac{i}{2}}$ and so also $R=\tau_{\frac{i}{2}}S$. As $S$ is a 
anti-homomorphism and $\tau_{\frac{i}{2}}$ a homomorphism, $R$ will be 
again a anti-homomorphism. As $S$ flips the comultiplication, and $\tau$ 
leaves it invariant, also $R$ will flip the comultiplication.
\snl
Finally, let us look at the behaviour of $R$ with respect to the 
involution. First, we have the general property
$$R(x^*)=S(\tau_{\frac{i}{2}}(x^*))
=S(\tau_{-\frac{i}{2}}(x)^*)
=S^{-1}(\tau_{-\frac{i}{2}}(x))^*$$
so that $R(x^*)=R^{-1}(x)^*$. In this particular case, and this is of 
course why $\tau_t$ has been defined that way, we have
$$\align R^2(a)&=S^2\tau_{i}(a)=S^2(a)=a\\
         R^2(b)&=S^2\tau_{i}(b)= e^{\frac{2\pi i}{n}}S^2(b)
	 =e^{\frac{2\pi i}{n}}a^{-1}ba=b.
\endalign$$
So $R^2=\iota$ and $R$ is involutive.  Together with the fact that 
$R(x^*)=R^{-1}(x)^*$, we see that $R$ is also a $^*$-map. This 
completes the result.
\einspr

The one-parameter group $(\tau_t)$ 
is called the {\it scaling group} while the map 
$R$ is usually called the {\it unitary antipode}.
Observe that we are using the convention $S=R\tau_{-\frac{i}{2}}$ (with 
the minus sign) to 
define the one-parameter group $\tau_t$. This is the convention used by 
Kustermans and Vaes in [K-V2] and it is different from the one used by 
Woronowicz in [W5].
\snl
We now calculate some more concrete formulas here as we will need them further 
in the paper.
\snl
We have $R(a)=a^{-1}$ and $R(a^*)={a^*}^{-1}$. We have also 
$$ R(b)=S\tau_{\frac{i}{2}}(b)=S(e^{\frac{\pi i}{n}}b)
=-e^{\frac{\pi i}{n}}a^{-1}b$$
and similarly
$$R(b^*)=S\tau_{\frac{i}{2}}(b^*)=S(\tau_{-\frac{i}{2}}(b)^*)
=S(e^{\frac{\pi i}{n}}b)=-e^{\frac{\pi i}{n}}{a^*}^{-1}b^*.$$
Observe further that
$$(e^{\frac{\pi i}{n}}a^{-1}b)^*=e^{-\frac{\pi i}{n}}b^*{a^*}^{-1}
=e^{-\frac{\pi i}{n}}e^{\frac{2\pi i}{n}}{a^*}^{-1}b^*$$
and indeed, we have $R(b)^*=R(b^*)$.
\snl
Observe that the unitary antipode is a complicated map. However, it can 
be seen as the composition of 3 maps:
$$\alignat 3  a &\to a^{-1} & \qquad\qquad  a &\to a 
& \qquad\qquad  a &\to a \\
b &\to b      & \qquad\qquad  b &\to -b 
&\qquad\qquad  b &\to e^{\frac{\pi i}{n}}a^{-1}b
\endalignat$$
(see section 3). The first one is a $^*$-anti-automorphism, while the 
second and the third ones are $^*$-automorphisms. The complicated part 
is put in the third one, but because now this is an automorphism and 
not a anti-automorphism, this becomes easier to treat.
\snl
We finish this section by saying briefly something about the underlying 
Hopf algebra for the other examples that treat, in a much more concise 
manner, in section 5.  One example there is very similar to the main 
example and has the same underlying Hopf $^*$-algebra, namely the one 
in 2.2, but now with a real deformation parameter. The Hopf algebra of 
2.1 will serve as a basis for the quantized $ax+b$-group and the 
quantized $E(2)$-group. However, different involutions are used. 	  
\nl\nl


\bf 3. The C$^*$-algebra and the comultiplication \rm
\nl
In the previous section we gave a complete description of the quantum 
$az+b$-group on the Hopf $^*$-algebra level.  In this section we will 
describe how this is lifted to the C$^*$-algebra level by Woronowicz in 
[W5].
\snl
The first step is the realization as operators on a Hilbert space 
of the generators $a$ and $b$ of the 
Hopf $^*$-algebra as given in proposition 2.3.  The next step will be to  
take appropriate functions of these operators.  This will be possible if 
we take nice representations. Then we get the associated C$^*$-algebra 
and von Neumann algebra. The comultiplication is defined by means of 
the implementation with the multiplicative unitary. Finally, at the end 
of the section, we study the antipode in this operator algebra 
framework.
\snl
As before, let $n$ be a non-zero natural number. We have seen  
(see the remark after proposition 2.3) that we 
should certainly take $n\geq 2$. In the introduction of the paper of 
Woronowicz ([W5]), where $N=2n$, it is mentioned that one 
should take $n\geq 3$, but it is claimed that the results are still 
valid for $n=2$. Therefore, we will just assume that $n\geq 2$.
\snl
We put $q = \exp \frac{\pi i}{n}$.  For convenience we will also use  
$q^{1/2} = \exp \frac{\pi i}{2n}$ and $q^{-1/2} = \exp(-\frac{\pi i}{2n})$ 
and similarly for other powers of $q$.
\snl
We will work with pairs $(a,b)$ of normal operators on a Hilbert space 
$\Cal H$.  In the following definition, we give the reformulation of 
the algebraic relations from the previous section in terms of these 
operators.

\inspr{3.1} Definition  \rm
Let $a = u|a|$ and $b = v|b|$ be the polar decompositions of the operators 
$a$ and $b$ respectively.  We will first of all assume that $|a|$ and $|b|$ 
are non-singular (i.e.\ have trivial kernel).  So, $u$ and $v$ are unitary.  
Furthermore we assume
\snl
\quad i) $u^{2n} = v^{2n} = 1$,

\quad ii) $|a|$ and $v$ commute; $|b|$ and $u$ commute,

\quad iii) $uv = qvu$,

\quad iv) $|a|^{it} |b| |a|^{-it} = e^{- \frac{\pi t}{n}} |b|$. 
\snl
We will say that the pair $(a,b)$ is an {\it admissible pair} of normal 
operators.
\einspr

These conditions are the same as the conditions 0.7, formulated in the 
paper of Woronowicz ([W5]), except for the fact that we also assume 
that $|b|$ is non-singular.
\snl
Before continuing, we collect some remarks about this definition.

\inspr{3.2} Remarks \rm 
i) It is quite natural to require $|a|$ to be non-singular because in the Hopf 
$^*$-algebra we have that the element $a$ is invertible.  It is not so obvious 
to take also $|b|$ non-singular.  Of course the requirement $v^{2n} = 1$ 
would not be possible and we would need to formulate this in a different way.  
It can be done (again see 0.7 in [W5]), 
but it turns out that there is no need for this.
\snl
ii) The requirements that $u^{2n} = 1$ and $v^{2n} = 1$ are equivalent 
with a restriction on the spectra of these operators.  Let
$$\Gamma = \{ q^k r \mid k = 0,1,\dots ,2n-1  \text{ and } r > 0 \}.$$
The requirement is that the spectrum of $a$ and $b$ belong to the 
closure $\overline \Gamma$ of $\Gamma$ (which is $\Gamma \cup \{ 0 \})$.   
\snl
iii) We have already that $u$ commutes with $|a|$ and that $v$ commutes 
with $|b|$.  We moreover assume that $|a|$ commutes with $v$ and that 
$|b|$ commutes with $u$.  So the pair $(|a|,|b|)$ commutes with the pair 
$(u,v)$.  This is the translation of the fact that, in the Hopf 
$^*$-algebra, the elements $a$ and $b$ are normal and that $a$ commutes 
with $b^*$ (and so $a^*$ commutes with $b$). 

\snl
iv) The conditions iii) and iv) translate the commutation rule 
$ab = q^2 ba$ which we had for the elements $a,b$ in the Hopf 
$^*$-algebra (indeed $q^2 = \lambda = e^{\frac{2 \pi i}{n}}$).  
The fact that it splits up as it does is a consequence of the condition 
$|\lambda| = 1$.  See section 5 for another case.

\snl
v) Two unitaries, $u$ and $v$, satisfying $u^{2n} = v^{2n} = 1$ and 
$uv = qvu$ where $q = e^{\frac{\pi i}{n}}$ always generate an algebra 
isomorphic with $M_n(\Bbb C)$, the $n \times n$ matrix algebra over the 
complex numbers.  In fact $u$ and $v$ induce representations of the 
group $\Bbb Z_{2n}$ and $v$ is considered as giving a representation of 
the dual group $\hat \Bbb Z_{2n}$ (which is of course $\Bbb Z_{2n}$) so 
that they combine to a  pair $(u,v)$ giving the Heisenberg representation of 
$\Bbb Z_{2n}$ and its dual $\hat \Bbb Z_{2n}$.  Here this means that 
$u^kv^\ell = q^{k \ell} v^\ell u^k$ where $(k,\ell) \rightarrow q^{k \ell}$ 
is precisely the bicharacter that defines the pairing between $\Bbb Z_{2n}$ 
and $\hat \Bbb Z_{2n}$. See the example A.4.i in the appendix.  
\snl
vi) We have a similar situation for the commutation rule between 
$|a|$ and $|b|$.  The form we have given is a strong form of the formal 
relation $|a||b| = q|b||a|$.  The Heisenberg form now is 
$|a|^{it} |b|^{is} = e^{-\frac{\pi ist}{n}} |b|^{is} |a|^{it}$ and 
again the map 
$(t,s)\to e^{-\frac{\pi ist}{n}}$ is a pairing between $\Bbb R$ and its dual 
$\hat \Bbb R$ (again $\Bbb R$).  Also here, such a representation is 
determined, up to multiplicity.  The von Neumann algebra generated by 
operators $|a|$ and $|b|$ like this is isomorphic with $\Cal B(\Cal K)$ 
for some separable Hilbert space $\Cal K$.  
If we take the C$^*$-algebra generated by elements 
$\int f(t) |a|^{it} \,dt$ and $\int g(s)|b|^{is} 
\,ds$ with $f$ and $g$  
continuous with compact support, we get a C$^*$-algebra isomorphic with 
the compact operators on $\Cal H$. We refer to the example A.4.ii in 
the appendix and the general remarks before the examples.
\einspr

It is good to have in mind that there is the group $\Gamma$ behind.  In 
fact $\Gamma = \Bbb Z_{2n} \times \Bbb R$ and so $\Gamma$ is again 
self-dual.  The representations that are given by $a$ and $b$ also form a 
Heisenberg type representation (cf.\ example A.4.iii).  
We will not take this point of view 
further in the paper although it should be mentioned that this would 
yield simpler formulas.  But this would be less familiar and therefore 
we will not use this.
\snl
Let us now, as an example, give the simple Heisenberg realisations of such 
operators.  These basic, irreducible representations will play a role 
further, when we obtain the right regular representation from the right 
Haar measure in section 4 (see proposition 4.15). 

\inspr{3.3} Proposition \rm
Consider the Hilbert space $\Bbb C^{2n}$ with a basis
$\{e_k \mid k = 0,1, \dots ,2n-1 \}$.  Define operators $m$ and $s$ by
$$\align
me_k & = q^k e_k \\
se_k & = e_{k+1}
\endalign
$$
where it is understood that $e_{2n} = e_0$.  Then $m$ and $s$ are unitary and 
$m^{2n} = s^{2n} = 1$ and $ms = qsm$.
\einspr

The result is well-known and the proof is trivial.  
It is also easy to  see that $m$ and $s$ generate $M_{2n}(\Bbb C)$.  
\snl
In what follows, we will agree that not only $e_{2n} = e_0$ but also 
that $e_{k+2n} = e_k$ for all $k\in \Bbb Z$ and hence $e_{-k} = 
e_{2n-k}$. This means that we consider the basis indexed over 
$\Bbb Z_{2n}$. 
\snl
Less trivial, but also well-known is the following.

\inspr{3.4} Proposition \rm  
Now take the Hilbert space $L^2(\Bbb R^+)$ where $\Bbb R^+$ is considered 
with the usual Lebesgue measure.  Define positive, self-adjoint, non-singular 
operators $a_0$ and $b_0$ on this Hilbert space by
$$\align
(b_0 f)(s) & = sf(s) \\
(a^{it}_0 f)(s) & = e^{- \frac{\pi t}{2n}} f(e^{- \frac{\pi t}{n}}s)
\endalign
$$
where $t \in \Bbb R$.

\einspr

Observe that the operator $a_0$ is defined as the analytic generator 
of a one-parameter group of unitaries.  We do not get the common form of 
the Heisenberg representation as this is defined on $L^2 (\Bbb R)$; but 
it is easy to construct the unitary  form one space to the other, 
relating the two forms.  It is the form we give that will appear 
naturally further (when we study the right Haar measure and the 
regular representations).

\snl
Of course, when we take the tensor product  
$\Cal H = L^2 (\Bbb R^+) \otimes \Bbb C^{2n}$
of these two Hilbert spaces, we get an admissible pair 
$(a,b)$ satisfying the assumptions by taking $a = a_0 \otimes m$ and 
$b = b_0 \otimes s$.
 In fact, it follows from the general theory (see 
the appendix), that any admissible pair is obtained from this 
irreducible one by taking the tensor product with $1$ on some other 
Hilbert space.

\nl
This is the {\it first step} in the process of passing from the Hopf 
$^*$-algebra to the C$^*$-algebra.  The {\it next step} is to consider 
appropriate functions.  And it is not only necessary  to take the correct 
functions, it is always important to have them in a form suitable for further 
reasoning and calculations.
\snl
Such functions are used in the following definition, preliminary to the 
introduction of the C$^*$-algebra.

\inspr{3.5} Proposition \rm  
Let $(a,b)$ be an admissible pair of normal operators on a Hilbert space 
$\Cal H$.   For $k, \ell = 0,1,\dots,2n-1$, let 
$f_{k,\ell}$ be a continuous complex function with compact support in 
$\Bbb R^+ \times \Bbb R$.  Assume 
that $f_{k,\ell} (0,t) = 0$ for all $t$ and all $\ell$ when $k \ne 0$.  
Elements of the form
$$\sum_{k,\ell} \left(\int f_{k,\ell} (|b|,t) |a|^{it} \,dt
\right) v^k u^\ell$$
form a non-degenerate $^\ast$-algebra $A_0$ of bounded operators on $\Cal H$.

\snl {\bf Proof:} 
This is essentially straightforward.  But as we will need an explicit 
formula for the adjoint and the product later, we will give it here anyway.  
So let $x,y \in A_0$ given by
$$\align
x = & \sum_{k,\ell} \left(\int f_{k,\ell} (|b|,t) |a|^{it} \,dt
\right)v^k u^\ell \\
y = & \sum_{k^\prime,\ell^\prime} \left(\int g_{k^\prime,\ell^\prime} 
(|b|,s)|a|^{is}\,ds\right)v^{k^\prime} u^{\ell^\prime}.
\endalign
$$
Then, using that  $|a|^{it}|b||a|^{it} = e^{-\frac{\pi t}{n}} |b|$ and 
$u^\ell v^{k^\prime} = q^{k^\prime \ell} v^{k^\prime} u^\ell$, we find for 
the product
$$xy = \sum_{k,\ell,k^\prime,\ell^\prime} \left(\iint f_{k,\ell} 
(|b|,t)g_{k^\prime,\ell^\prime} (e^{-\frac{\pi t}{n}}|b|,s) 
|a|^{i(t+s)} \,dt \,ds\right) q^{k^\prime \ell} v^{k + k^\prime} u^{\ell + 
\ell^\prime}.$$
This is again an element of the same type.  Observe that, when 
$k + k^\prime \ne 0$ then either $k \ne 0$ or $k^\prime \ne 0$ and so 
evaluation of the first variable in $0$ will give $0$ as this is the case 
for $f_{k,\ell}$ or for $g_{k^\prime,\ell^\prime}$.
\snl
Similarly, using the same relations,
$$x^\ast = \sum_{k,\ell} \left(\int \overline{f_{k,\ell}} (e^{\frac{\pi t}{n}} 
|b|, t) |a|^{-it} \,dt\right) q^{k\ell} v^{-k} u^{-\ell}$$
which is again in the algebra.  This proves the proposition.
\einspr     

Of course, we have a $^*$-representation of the crossed product of the 
C$^*$-algebra $C_0(\Bbb R_+)$ by the action $\alpha$ of $\Bbb R$ given 
by $(\alpha_{t}f)(p) = f(e^{-\frac{\pi t}{n}}p)$.  
And we have the tensor product with $M_{2n}$, given by the operators $u$ 
and $v$.  However, we take a certain $^\ast$-subalgebra because of our 
condition $f_{k,\ell}(0,t) = 0$ for all $t$ and all $\ell$ when $k \ne 0$.  
We will explain later where this condition comes from and why it is important
(see the remark after proposition 3.8).

\inspr{3.6} Definition \rm  
Let $A$ be the C$^*$-algebra on $\Cal H$ obtained by taking the norm closure 
of the $^\ast$-algebra $A_0$ in the previous proposition.
\einspr

Observe that this C$^*$-algebra is independent of the choice of the pair 
$(a,b)$.  In fact, any such pair will give a faithful representation of 
this C$^*$-algebra.  So we can think of the C$^*$-algebra as an abstract 
C$^*$-algebra. It is a certain subalgebra of the tensor product of the crossed 
product $(C_0(\Bbb R^+) \times_\alpha \Bbb R) \otimes M_{2n}$. That it is 
a proper subalgebra is coming from the fact that a function of the 
factors  $|b|$ and $v$ in the polar decomposition of $b$ must really 
come from a function on $b$. 
Any admissible pair $(a,b)$ will give a faithful 
non-degenerate $^*$-representation of this C$^*$-algebra.

\snl
We will also use the von Neumann algebra generated by this 
C$^*$-algebra. We will denote it by $M$. 
It contains the multiplier algebra $M(A)$ of $A$. 
It is isomorphic with the 
tensor product $\Cal B(\Cal K)\ot M_{2n}(\Bbb C)$ where $\Cal K$ is any 
separable Hilbert space. This von Neumann algebra is determined up to a 
multiplicity, depending on the choice of the pair $(a,b)$. Also observe 
that, on this von Neumann algebra level, the restriction that we have 
on the functions $f_{k,\ell}$ in proposition 3.5 is no longer 
important. In some sense, it would be much more easy to only look at the von 
Neumann algebra framework here. On the other hand, from a theoretical 
point of view, the C$^*$-algebra approach is more natural. Moreover, 
this is the framework used by Woronowicz in [W5]. 
\nl
Now, we turn our attention to the {\it comultiplication}.  The direct way would 
be to start with an admissible pair $(a,b)$ of normal operators on $\Cal H$ 
and define a new pair $(\tilde a, \tilde b)$  on 
$\Cal H \otimes \Cal H$ by 
$$\align
\tilde a = & a \otimes a \\
\tilde b = & a \otimes b + b \otimes 1
\endalign
$$
where in fact 
we take the closure of this sum.  And indeed, it is true that this 
gives again a pair of normal operators with the right spectral conditions 
and commutation rules.
\snl
This track however, though natural at the first thought, is not very easy.  
It turns out that it is easier to construct a good multiplicative 
unitary.  Such a multiplicative unitary can be found by first looking at 
unitary operators $W$ on $\Cal H \otimes \Cal H$ such that
$$\align
\tilde b & = W(b \otimes 1) W^\ast \\
\tilde a & = W(a \otimes 1) W^\ast
\endalign
$$
In [VD4], we give a technique to construct such a $W$ from a dual pairing.  
And we apply this technique, as far as possible, to the examples, studied 
in this paper.

\snl
Here, we simply will describe the {\it multiplicative unitary}, as it was 
discovered in the paper by Woronowicz [W5].  However, we will not give 
the complete definition; we refer to [W5] for this.  Instead, we will 
concentrate on its properties and state them when we need them.  Also, 
we will use a slightly different (and somewhat simpler) multiplicative 
unitary which turns out to be sufficient for our purpose.
\snl
The first ingredient to construct this multiplicative unitary is a form of 
the quantum exponential function.  It is a continuous function $F$ defined 
on the group $\Gamma = \{q^k r \mid k = 0,1,\dots,2n-1 \text{ and } r > 0\}$, 
with values  in the unit circle of $\Bbb C$.  See equation 1.5 in   [W5] 
where the function is called $F_N$ ($N$ is $2n$ here).  The main properties 
of this function are collected in proposition 1.1 of [W5].  The second 
ingredient is a bicharacter $\chi$ on $\Gamma \times \Gamma$.  It is defined 
as
$$\chi (\gamma,\gamma^\prime) = q^{kk^\prime} e^{\frac{n}{\pi i} 
(\log r)(\log r^\prime)}$$
where $\gamma = q^k r$ and $\gamma^\prime = q^{\ell^\prime}r^\prime$.  
Again see formula 1.1 in [W5].  This bicharacter is essentially a pairing 
realizing the self-duality of $\Gamma$ which we mentioned already.
\snl
Then we are ready to recall the following definition of Woronowicz 
(theorem 3.1 in [W5]).

\inspr{3.7} Definition \rm  
Consider two admissible pairs $(a,b)$ and $(\hat a, \hat b)$ of 
normal operators on 
Hilbert spaces $\Cal H$ and $\Cal K$ respectively.  
Then we define a unitary $W$ on $\Cal K \otimes \Cal H$ 
by
$$W = F(\hat b \otimes b)\chi (\hat a \otimes 1, 1 \otimes a).$$
\einspr

Recall that the spectra of $a,\hat a,b$ and $\hat b$ are contained in 
$\overline\Gamma = \Gamma \cup \{0 \}$.  Then the same is true for 
$\hat b \otimes b$, $\hat a \otimes 1$ and $1 \otimes a$.  As these 
operators are non-singular, $0$ is not in the point spectrum.  Therefore, 
we can apply the function $\chi$ which is only defined on 
$\Gamma \times \Gamma$.  There is no problem with applying $F$ because 
this is defined on $\overline\Gamma$.

\snl
As $F$ and $\chi$ map into the unit circle, $W$ is a unitary.
\snl
We will use this unitary for various choices of the pair $(\hat a, \hat b)$.  
In general, we can already state the following result.

\inspr{3.8} Proposition \rm  
Let $(a,b)$ and $(\hat a, \hat b)$ and $W$ be as in definition 3.7.  
Let $A$ be the C$^*$-algebra as defined in 3.6 for the pair $(a,b)$.  
Then $A$ is the norm closure of the set 
$\{ (\omega \otimes \iota) W \mid \omega \in \Cal B(\Cal K)_\ast \}$.
\einspr

We use $\Cal B(\Cal K)_\ast$ for the predual of $\Cal B(\Cal K)$, the space of 
normal linear functionals on $\Cal B(\Cal K)$.   
\snl  
This result as such is not stated in [W5] but it can be deduced from 
properties that are proven by Woronowicz.  We refer to section 6 of [W5].  
There, a special choice of $(\hat a, \hat b)$ is considered, but it is 
clear that this is not important for the result above.  Then formula 6.5 
together with formula 4.2 of [W5] will give the result. Observe that 
the C$^*$-algebra $A$ is the crossed product of the C$^*$-algebra 
$C_0(\overline \Gamma)$ by the action $\sigma$ of the group $\Gamma$ 
defined by $(\sigma_\gamma f)(\gamma')=f(\gamma\gamma')$ whenever 
$\gamma,\gamma'\in \Gamma$. This clarifies the restriction on the 
functions $f_{k,\ell}$ that we have in definition 3.5.
\snl
We have that $W\in M(\Cal B_0(\Cal H)\ot A)$ where we use 
$\Cal B_0(\Cal H)$ to denote the C$^*$-algebra of compact operators on 
the Hilbert space $\Cal H$. Of course, also $W\in \Cal B (\Cal H)\ot M$ 
where the von Neumann algebra tensor product is considered.
\snl
Because of the symmetry, we have similar properties for the other leg. 
In fact, when we use $\hat A$ and $\hat M$ for the C$^*$-algebra and 
von Neumann algebras associated with the pair $(\hat a,\hat b)$, we get 
$W\in M(\hat A\ot A)$ and also $W\in \hat M\ot M$ (with the von Neumann 
algebra tensor product).  
\snl
There is another property of this unitary concerning the antipode.  We 
will come back to this later (see 3.11 and 3.12).
\snl
Observe that in this special example, we have that $\hat A$ is 
isomorphic with $A$ and that the same is true for $\hat M$ and $M$. This is 
typical for a self-dual example such as the one studied here. Also two 
other examples, studied in section 5 are self-dual. However, the 
quantum $E(2)$, that we only will discuss very briefly in section 5, is 
not self-dual and therefore does not have this property. This 
observation is important and we will come back to it in the last 
section where we draw some conclusions.
\nl
The next step is the definition of the comultiplication.  We first need the 
following lemma.

\inspr{3.9} Lemma \rm  
Let $(a,b)$ an admissible pair of normal operators on a Hilbert 
space $\Cal H$.  Let $\hat a = b^{-1}$ and let $\hat b$ be 
the closure of $ab^{-1}$.  
Then $(\hat a, \hat b)$ is also an admissible pair of normal operators.

\snl \bf Proof: \rm
The polar decomposition of $\hat a$ is clearly $v^\ast |b|^{-1}$ 
(where $b = v|b|$ is the polar decomposition of $b$ as before).  When 
$a = u|a|$ is the polar decomposition of $a$, then we write
$$ab^{-1} = uv^\ast |a||b|^{-1} = (q^{-1/2} uv^\ast)(q^{1/2} |a||b|^{-1}).$$
We know that $q^{1/2} |a||b|^{-1}$ is self-adjoint and positive (cf. 
proposition A.5 for a similar result).  
And $q^{-1/2} uv^\ast$ is unitary.  So this is the polar decomposition of 
$ab^{-1}$. Now it is straightforward to verify the assumptions.
\einspr

Observe that we have used (as we will do further in the paper) 
$ab^{-1}$ and similarly $|a||b|^{-1}$ to denote the closures of these 
operators.
\snl
Finally,  we are ready to give the definition of the comultiplication. 
        	  
\iinspr{3.10} Proposition \rm
Let $(a,b)$ and $(\hat a,\hat b)$
 be as in the previous lemma and let 
$W$ be the unitary as defined in 3.7. Then $W$ is a multiplicative 
unitary and $\Phi$ defined on $A$ by $\Phi(x)=W(x\ot 1)W^*$ is a 
comultiplication on $A$.
\einspr
Again, this is sligthly different from what is found in theorem 5.1 of 
[W5]. There, $\hat a$ is defined as $sb^{-1}$ where $s$ is a 
non-singular positive self-adjoint operator that strongly commutes with 
$a$ and $b$. This is needed to prove that $W$ is manageable. We will 
not need this property here. We can therefore take $s=1$ which amounts 
to applying a C$^*$-homomorphism and will respect the multiplicative 
property. See also formula 0.11 in [W5].
\snl
Remark that it follows from the pentagon equation that $\Phi(A)\in 
M(A\ot A)$ because $W\in M(\Cal B_0(\Cal H)\ot A)$. Moreover, because 
elements of the form $(\omega\ot\iota)W$ with $\omega\in \Cal B(\Cal 
H)_*$ belong to $A$ and form a dense subspace of it, slices of 
$\Phi(A)$ with functionals $\omega \in A^*$ lie in $A$ and both left 
and right slices give dense subspaces of $A$. This property is 
important to have a locally compact quantum group in the next section 
(see theorem 4.11).
\snl
Further, remark that the comultiplication is also defined on the von 
Neumann algebra $M$ by the same formula as in the proposition and that 
$\Phi(M)$ lies in the von Neumann algebra tensor product $M\ot M$. 
\snl
We have made a very special choice in proposition 3.10. We have 
expressed the generators of $\hat A$ in terms of the generators of $A$. 
Now, both C$^*$-algebras are acting on the same Hilbert space. They are 
isomorphic as we mentioned before, but they are not equal. On the other 
hand, the von Neumann algebras $\hat M$ and $M$ now coincide, in fact 
with all of $\Cal B(\Cal H)$. Again this is typical for the 
self-duality of this example. However, there is more. The fact that 
$\hat M=M=\Cal B(\Cal H)$ means that the multiplicative unitary $W$ in 
proposition 3.10 can not be the regular representation because then one 
has $\hat M \cap M = \Bbb C 1$. The choice that Woronowicz has made 
(using a non-trivial operator $s$ - see above) gives an intermediate 
situation. See also the remarks that we will make further when we 
construct the regular representation. A similar situation occurs with 
the two other examples in section 5. The situation is again different 
with the quantum $E(2)$ as this is not self-dual. We will come back to 
this remark in section 6. 
\nl
Finally, we look at the {\it antipode}. From the general theory, we know that 
the antipode $S$ is, roughly speaking, characterized by the equation 
$(\iota\ot S)W=W^*$. In fact, the following result is more or less 
obvious (and standard).
\snl
Before we formulate this result, we have to make an {\it important 
remark}. 
There is, in the general theory, a choice to make for the 
comultiplication on the dual. In the algebraic framework (e.g.\ [VD1], 
[VD3] and [K-VD]), it is quite common to define the comultiplication 
dual to the multiplication. In the C$^*$-framework (e.g.\ [K-V2]), 
it is common to take the opposite comultiplication on the dual. This 
choice is also made by Woronowicz. We however have made the first choice 
and just taken the comultiplication dual to the multiplication as in the 
algebraic setting. The difference between the two choices has e.g.\ a 
consequence for the antipode $\hat S$. One needs to take the inverse 
for the other choice.

\iinspr{3.11} Lemma \rm
Consider the unitary $W$  as defined in 3.7 for any two pairs $(a,b)$ 
and $(\hat a,\hat b)$. Then, there are linear, densely defined maps $S$ 
from $A$ to $A$  and $\hat S$ from $\hat A$ to $\hat A$ given by the 
formulas
$$\align
S((\omega\ot\iota)W) &=(\omega\ot\iota)W^*\\
\hat S((\iota\ot\omega)W) &=(\iota\ot\omega)W^*
\endalign$$ 
where $\omega$ runs through the predual $\Cal B(\Cal H)_*$ of 
$\Cal B(\Cal H)$ in the first formula and through $\Cal B(\Cal K)_*$ in 
the second case.

\snl\bf Proof: \rm
The argument is simple. One just has to argue that 
$(\omega\ot\iota)W^*=0$ when $(\omega\ot\iota)W=0$. Now, 
$(\omega\ot\iota)W^*=((\overline\omega\ot\iota)W)^*$ where 
$\overline\omega(x)=\omega(x^*)^-$ for $x\in \Cal B(\Cal H)$. 
Therefore, the result follows from the fact that the norm closure of 
the space $\{(\iota\ot\omega)W=0 \mid \omega\in \Cal B(\Cal H)_* \}$
is self-adjoint. Similarly on the other side.
\einspr

This general result is not very useful. One would like to have more 
information about $S$. Now, we have seen in the algebraic setting in 
the previous section how $S$ is defined on the elements $a$ and $b$ and 
that it has a polar decomposition. Therefore, the following result is 
no surprise.

\iinspr{3.12} Proposition \rm
Let $W$ be as in 3.7 and let $A$ be the C$^*$-algebra associated with 
$(a,b)$ as in 3.6. There exists a strongly continuous one-parameter 
group of $^*$-automorphisms $(\tau_t)_{t\in \Bbb R}$ of $A$ 
characterized by the action on the generators, denoted and given by 
$\tau_t(b)=e^{\frac{2\pi t}{n}}b$ and $\tau_t(a)=a$. There exists an 
involutive $^*$-anti-automorphism $R$ of $A$, also characterized by the 
action on the generators given by $R(a)=a^{-1}$ 
and $R(b)=-qa^{-1}b$. We have that $R$ commutes with $\tau$. Also
$$\align
\Phi(\tau_t(x))&=(\tau_t\ot \tau_t)\Phi(x)\\
\Phi(R(x))&=\sigma(R\ot R)\Phi(x)
\endalign$$
for any $x\in A$ (where $\sigma$ is the flip). Moreover, let $S$ be 
defined by $R\tau_{-\frac{i}{2}}$ where $\tau_{-\frac{i}{2}}$ is the 
analytic extension of $(\tau_t)$ to the point $-\frac{i}{2}$ (as an 
unbounded map of course). Then, $(\omega\ot\iota)W$ belongs to the 
domain of $S$ and $S((\omega\ot\iota)W)=(\omega\ot\iota)W^*$. In fact, 
we get a core in the domain.
\einspr

As before, $R$ is called the unitary antipode and $\tau$ the scaling 
group. Recall that we use a different convention for the polar decomposition 
of the antipode than Woronowicz (see also the remark at the end of the previous 
section).
\snl
By symmetry, similar data exist for the algebra $\hat A$ associated 
with the pair $(\hat a,\hat b)$. The formulas are the same.
\snl
As we mentioned already, the result is not unexpected. However, the 
proof is far from trivial and it can be found in [W5]. 
It mainly uses very detailed analysis of the 
function $F$ and its Fourier transform.
\snl
We will not give the proof here of course, but let us nevertheless show 
how to obtain all but the last property in a standard way.
\snl 
A first step is the following lemma.

\iinspr{3.13} Lemma \rm
The $^*$-automorphisms $\tau_t$ are given by 
$\tau_t(x)=|a|^{-2it}x|a|^{2it}$.
   
\snl\bf Proof: \rm 
Observe that $|a|^{-2it}a|a|^{2it}=a$ and 
$|a|^{-2it}b|a|^{2it}=e^{\frac{2\pi t}{n}}b$ by the assumptions in 3.1. 
Then, it is immediately clear that this gives a strongly continuous 
one-parameter group of $^*$-automorphisms of $A$.
\einspr

Now we will have $\Phi(|a|^{it})=|a|^{it}\ot |a|^{it}$ (recall that the 
comultiplication is also defined on the multiplier algebra $M(A)$ and even on 
the von Neumann algebra $M$). Then it follows 
that $\Phi(\tau_t(x))=(\tau_t\ot\tau_t)\Phi(x)$.
\snl
The situation with $R$ is more complicated. But it is possible to 
obtain $R$ as a composition of three maps which are easier to 
understand. Let us look again at the algebraic situation (cf.\ the end 
of section 2). The map $R$ is a 
anti-automorphism characterized by $R(a)=a^{-1}$ and $R(b)=-qa^{-1}b$. 
This is the composition of the following three maps :
\snl
\quad i) $R_1$ mapping $a$ to $a^{-1}$ and $b$ to $b$,

\quad ii) $R_2$ mapping $a$ to $a$ and $b$ to $-b$,

\quad iii) $R_3$ mapping $a$ to $a$ and $b$ to $qa^{-1}b$.
\snl
All these maps commute and so the order is not important. We have that 
$R_1$ is a 
$^*$-anti-automorphism while the two other maps are $^*$-automorphisms.
\snl
All of these maps are implemented when the elements are realized as 
operators in 3.3 and 3.4. Then the first mapping is realized as $x\to 
Jx^*J$ where $J$ is a anti-linear map, defined as $J_0\ot J_1$ where
$$\align (J_0f)(s)&=\overline f(s) \qquad\qquad f\in L^2(\Bbb R^+)\\
         J_1e_k&=e_{-k} \qquad\qquad k=0,1,\dots,2n-1
\endalign$$
where, as before, we consider the last basis of $\Bbb C^{2n}$ as indexed 
over $\Bbb Z_{2n}$.
\snl
The second map is realized by the unitary $1\ot v_0$ where $v_0$ is 
defined on $\Bbb C^{2n}$ by $v_0e_k=(-1)^ke_k$. This $^*$-automorphism 
is inner.
\snl
Again, the last map is the more complicated one. It is again inner and 
the result follows from the following lemmas.

\iinspr{3.14} Lemma \rm
Let $h=\log |a|$ and $u_0=\exp\frac{in}{2\pi}h^2$. Then 
$u_0|b|u_0^*=q^{\frac12}|a|^{-1}|b|$.

\snl\bf Proof: \rm
We have to verify that 
$$u_0|b|^{it}u_0^*=e^{-\frac{\pi i t^2}{2n}}|a|^{-it}|b|^{it}$$
as the right hand side is the $it$-th power of the positive 
non-singular operator $q^{\frac12}|a|^{-1}|b|$ (see proposition A5 in 
the appendix). This equation can be rewritten as 
$$|b|^{it}u_0^*|b|^{-it}=e^{-\frac{\pi i t^2}{2n}}|a|^{-it}u_0^*.$$
But $|b|^{it}|a||b|^{-it}=e^{\frac{\pi t}{n}}|a|$ and so
$|b|^{it}h|b|^{-it}=h+\frac{\pi t}{n}$ and so 
$$\align |b|^{it}(-\frac{in}{2\pi}h^2)|b|^{-it} 
&=-\frac{in}{2\pi}(h+\frac{\pi t}{n})^2\\
&=-\frac{in}{2\pi}h^2-ith-\frac{\pi i t^2}{2n}
\endalign$$
and taking the exponential gives the required equation.
\einspr

This will take care of the positive part. The following will work for 
the unitary part.

\iinspr{3.15} Lemma \rm 
There is a unitary $v_1$ on $\Bbb C^{2n}$ such that $v_1$ commutes with 
$u$ and $v_1vv_1^*=q^\frac12 u^*v$.

\snl\bf Proof: \rm
Consider again the representation given in 3.3. Define $v_1e_k=c_ke_k$ 
where $c_k=q^{-\frac12 k^2}=e^{-\frac{\pi i k^2}{2n}}$. 
Then clearly $v_1$ commutes with $u$. 
Moreover 
$$v_1vv_1^*e_k=\overline{c_k} v_1ve_k=\overline{c_k} v_1e_{k+1}=
\overline{c_k} c_{k+1}e_{k+1}.$$
Now, one verifies that $\overline c_kc_{k+1}=q^{-k-\frac12}$ so that 
indeed $v_1vv_1^*=q^{\frac12}u^*v$.
\einspr

Taking the two results together, we see that $u_0\ot v_1$ implements 
the $^*$-automorphism $R_3$. This can also serve as a way to prove the 
existence of $R$ and to show its basic properties. However, to show that $R$ 
flips the comultiplication does not seem to be easy.
\nl\nl


\bf 4. The right Haar measure and the regular representation \rm
\nl
Consider the C$^*$-algebra $A$ and the comultiplication $\Phi$ as 
described in the previous section. We will construct a faithful, lower 
semi-continuous, densely defined KMS-weight $\psi$ on $A$ and prove that 
it is right invariant. We will also prove the existence of such a left 
invariant weight and we will show that the pair $(A,\Phi)$ becomes a 
locally compact quantum group in the sense of Kustermans and Vaes [K-V2]. 
The Haar weights are not invariant, but only relatively invariant with 
respect to the scaling group. Finally, we will construct the right regular 
representation and discuss the relation of this multiplicative unitary 
with the original one as given in the previous section.
\snl
We will freely use the notations of the previous section. In 
particular, $(a,b)$ is an admissible pair of normal operators on a 
Hilbert space $\Cal H$ and $A$ is the associated C$^*$-algebra as 
defined in 3.6. Similarly, $M$ is the associated von Neumann algebra, 
defined as the weak closure of $A$.
\snl
We begin with the construction of a positive linear functional $\psi$ 
on the $^*$-subalgebra $A_0$ of $A$ as defined in proposition 3.5.

\inspr{4.1} Proposition \rm
Define a linear functional $\psi$ on $A_0$ by 
$$\psi(x)=\int f_{0,0}(r,0)r\,dr$$
when 
$$x=\sum_{k,\ell}\left(\int 
f_{k,\ell}(|b|,t)|a|^{it}\,dt\right)v^ku^\ell$$
with the $f_{k,\ell}$ as in proposition 3.5.
Then $\psi$ is faithful and positive. With $x$ as before, we get
$$\psi(x^*x)= \sum_{k,l}\iint 
|f_{k,\ell}(r,t)|^2 e^{-\frac{2\pi t}{n}}r\,dr\,dt.$$
\einspr

Observe that the variable $r$ lies in $\Bbb R^+$ and that we use the 
Lebesgue measure on $\Bbb R^+$. The other variable $t$ lies in $\Bbb R$ 
and we integrate it with respect to the Lebesgue measure on $\Bbb R$. 
We will use these conventions everywhere in this section. As before, 
the indices $k$ and $\ell$ run over $\Bbb Z_{2n}$.

\inspr{} Proof: \rm
Let $x$ be as above. Using formulas given in the proof of 
proposition 3.5, we find
$$x^*x=\sum_{k,l,k',\ell'}\left(\iint 
\overline{f_{k,\ell}}(e^{\frac{\pi t}{n}}|b|,t) 
f_{k',\ell'}(e^{\frac{\pi t}{n}}|b|,s)
|a|^{i(s-t)}\,ds\,dt\right)
q^{-\ell(k'-k)}v^{k'-k}u^{\ell'-\ell}.$$
So, using the formula for $\psi$ in the formulation, we get
$$\align
\psi(x^*x)&=\sum_{k,l}\left(\iint \overline{f_{k,\ell}}(e^{\frac{\pi t}{n}}r,t) 
            f_{k,\ell}(e^{\frac{\pi t}{n}}r,t) r\,dr\,dt\right)\\
          &=\sum_{k,l}\iint 
|f_{k,\ell}(r,t)|^2 e^{-\frac{2\pi t}{n}}r\,dr\,dt.
\endalign$$
The positivity and the faithfulness follow immediately from this 
formula. This proves the proposition.
\einspr

Next, we consider the associated GNS-representation of $A_0$.

\inspr{4.2} Proposition \rm
Define a Hilbert space by
$$\Cal H_\psi=L^2(\Bbb R^+)\ot L^2(\Bbb R) \ot \Bbb C^{2n} \ot \Bbb 
C^{2n}.$$
Also define a linear map 
$\eta_\psi:A_0\to \Cal H_\psi$ by
$\eta_\psi(x)=\xi$ where $x$ is as before and
$$\xi(r,t)=
\sum_{k,\ell}e^{-\frac{\pi t}{n}}r^{\frac12}f_{k,\ell}(r,t) e_k\ot 
e_\ell.$$
Here we identify the Hilbert space $\Cal H_\psi$ with 
$L^2(\Bbb R^+\times \Bbb R, \Bbb C^{2n} \times \Bbb C^{2n})$
and use an orthonormal basis $(e_k)$ in $\Bbb C^{2n}$ as before. 
Then we have $\psi(x^*x)=\langle\eta_\psi(x),\eta_\psi(x)\rangle$ for 
all $x\in A_0$. Left multiplication gives a $^*$-representation 
$\pi_\psi$ of 
$A_0$ by means of bounded operators. This representation is 
characterized by the action of the generators. Denoting these also by 
$\pi_\psi(a)$ and $\pi_\psi(b)$, we have
$$\align
\pi_\psi(a)&=a_0 \ot a_1 \ot m \ot s \\
\pi_\psi(b)&=b_0 \ot 1 \ot s \ot 1 
\endalign$$
where $a_0$ and $b_0$ are as in proposition 3.4 and $m$ and $s$ are as in 
proposition 3.3, and where $a_1$ is defined on $L^2(\Bbb R)$ by 
the formula $(a_1^{is}\xi)(t)=\xi(t-s)$.

\einspr

First recall that, by definition, 
$\pi_\psi(x)\eta_\psi(y)=\eta_\psi(xy)$ whenever $x,y\in A_0$. Then
the proof of the proposition is very straightforward and the calculations 
are easy using the formula for $\eta_\psi$ that we have in the 
formulation of the proposition.
\snl
The positive non-singular self-adjoint operator $a_1$ will satisfy 
$(a_1\xi)(t)=\xi(t+i)$ for $\xi$ in the appropriate domain, extended 
analytically to an horizontal strip. 
\snl
Observe that in this proposition, we have the $^*$-representation 
$\pi_\psi$ of the algebra $A_0$ and that we characterize it by saying 
what it does on the generators. It is not so hard to see how this 
should be done. In fact, we have done something like this already in 
the formulation of proposition 3.12 in the previous section, when we 
introduced the polar decomposition of the antipode. 
We will do similar things on other occasions further in 
this section. In particular, as before, we will also characterize $^*$-automorphisms of $A$ by 
saying what they do on the generators. We will always use the same 
symbol. We are aware of the fact that this has to be done with some 
care, but it is clear that it does not cause any difficulty. 
To do it completely rigourously would 
just involve more arguments (and more complicated formulas) 
and we are afraid they would not greatly 
clarify what we are doing. We refer to the remark about this point of 
view made in the introduction.
\snl
The next step is of course the construction of a left Hilbert algebra, 
coming from the above GNS-representation of $A_0$.

\inspr{4.3} Proposition \rm The subspace $\eta_\psi(A_0)$ of $\Cal 
H_\psi$ is a left Hilbert algebra when it is given the $^*$-algebra 
structure inherited from $A_0$. When we use $T$ for the closure of the map 
$\eta_\psi(x)\mapsto\eta_\psi(x^*)$ where $x\in A_0$, then the polar 
decomposition $J|T|$ of $T$ is given as follows.
We can write $J=J_{01}\ot J_{23}$ where $J_{01}$ acts on $L^2(\Bbb 
R^+\times \Bbb R)$ and $J_{23}$ on $\Bbb C^{2n}\ot\Bbb C^{2n}$ as
$$\align(J_{01}\xi)(r,t)&=e^{-\frac{\pi t}{2n}} \overline \xi(e^{-\frac{\pi 
t}{n}}r,-t)\\
J_{23}(e_k\ot e_\ell)&=q^{k\ell}e_{-k}\ot e_{-\ell}.
\endalign$$
Furthermore, we can write $|T|=1\ot c_1^{-1} \ot 1\ot 1$ where $c_1$ 
is defined on $L^2(\Bbb R)$ by  
$$(c_1^{is}\xi)(t)= e^{-\frac{\pi ist}{n}}\xi(t).$$

\snl \bf Proof: \rm
Most of the axioms of a left Hilbert algebra (cf.\ [St]) 
are more or less obvious. 
They come essentially free with the above construction. We just need to 
show that the map $\eta_\psi(x) \mapsto \eta_\psi(x^*)$ is closable.
So, let $x$ be as before and $\xi=\eta_\psi(x)$ and denote $\xi^\sharp=T\xi=
\eta_\psi(x^*)$. 
\snl
Recall the convention for the orthonormal basis. We consider it as 
indexed over the group $\Bbb Z_{2n}$ as in section 3. Also in what 
follows, we will use the notation 
$$\xi=\sum\xi_{k,\ell}\ot e_k\ot e_\ell$$
for elements in $\Cal H_\psi$. Now $\xi_{k,\ell}\in L^2(\Bbb R^+ \times 
\Bbb R)$.
\snl
Then, from the formula for $x^*$ in the proof of proposition 3.5 we get
$$x^*=\sum_{k,\ell}\left(\int 
\overline{f_{k,\ell}} (e^{-\frac{\pi t}{n}}|b|,-t)|a|^{it}
\,dt\right)q^{k\ell}v^{-k}u^{-\ell}$$
and so
$$\align
\xi^\sharp(r,t)&=
\sum_{k,\ell} e^{-\frac{\pi t}{n}}r^{\frac12} \overline {f_{k,\ell}}
(e^{-\frac{\pi t}{n}}r,-t)q^{k\ell}e_{-k}\ot e_{-\ell}\\
&=\sum_{k,\ell} e^{-\frac{\pi t}{n}}r^{\frac12} e^{-\frac{\pi t}{n}}
e^{\frac{\pi t}{2n}}r^{-\frac12}\overline {\xi_{k,\ell}}
(e^{-\frac{\pi t}{n}}r,-t)q^{k\ell}e_{-k}\ot e_{-\ell}\\
&=\sum_{k,\ell}e^{-\frac{3\pi t}{2n}}\overline {\xi_{k,\ell}}
(e^{-\frac{\pi t}{n}}r,-t)q^{k\ell}e_{-k}\ot e_{-\ell}.
\endalign$$
Then we see that this operator has a closure $T$ and that the operators $J$ 
and $|T|$ in the polar decomposition $T=J|T|$ of $T$ are as in the 
formulation of the proposition.
\einspr

We use $T$ for this map and $J|T|$ for its polar decomposition because the 
more common symbols $S$ and $S=J\Delta^\frac12$ are used for other 
things in this paper.
\snl
Finally, we use the general theory for constructing our weight from 
this left Hilbert algebra (see e.g.\ [St]). This results in the following.

\inspr{4.4} Theorem \rm 
There is a faithful, lower semi-continuous densely defined KMS-weight 
$\psi$ on $A$ extending the linear functional on $A_0$ as defined in 
proposition 4.1. It is KMS with respect to 
the automorphism group $\sigma$ on A, 
characterized by the images of the generators, when using the notations
$\sigma_t(a)$ and $\sigma_t(b)$, given by 
$$\align
\sigma_t(a)&=e^{\frac{2\pi t}{n}}a\\
\sigma_t(b)&=b.
\endalign$$

\snl \bf Proof: \rm
Consider the faithful, normal semi-finite weight $\widetilde\psi$ on 
$\pi_\psi(A_0)''$ as obtained from the left Hilbert algebra, using the 
general procedure (see e.g.\ [St]). The restriction of this weight to 
the C$^*$-algebra $\pi_\psi(A_0)^-$ (where the norm closure is 
considered), is a faithful, lower semi-continuous weight. The 
representation $\pi_\psi$ extends from $A_0$ to an isomorphism of $A$ 
to this C$^*$-algebra. We can use this isomorphism to get the faithful 
lower semi-continous weight $\psi$ on $A$. It is clear from the 
construction that this weight extends the linear functional on $A_0$ as 
defined in proposition 4.1. In particular, the weight will be densely 
defined.
\snl
Now, a straightforward calculation shows that
$$|T|^{2it}\pi_\psi(a)|T|^{-2it}=e^{\frac{2\pi t}{n}}\pi_\psi(a).$$ 
The  reason is essentially that also the pair $(a_1, c_1)$ satisfies a 
Heisenberg type relation (see appendix, example A.4.iv) 
$$c_1^{-it}a_1c_1^{it}=e^{\frac{\pi t}{n}}a_1.$$
Similarly $|T|^{2it}\pi_\psi(b)|T|^{-2it}=\pi_\psi(b)$.
Now it is known that the weight $\widetilde\psi$ on the von Neumann algebra 
$\pi_\psi(A_0)''$ is KMS with respect to the modular automorphism 
group, implemented by the unitaries $|T|^{2it}$. Then it is clear that 
$\sigma$, as defined in the proposition, gives a strongly continuous 
one-parameter group of $^*$-automorphisms of the C$^*$-algebra $A$ and 
that this is essentially the restriction of the modular automorphism 
group above. Then it follows that we have a KMS-weight for this 
automorphism group $\sigma$. This completes the proof.
\einspr

It is important for what follows to keep in mind that the weight $\psi$ 
is obtained from restricting the faithful normal semi-finite weight 
$\widetilde\psi$ on the von Neumann algebra $\pi_\psi(A_0)''$. We know from the 
observations in the previous section that this von Neumann algebra is 
isomorphic with the von Neumann algebra $M$, generated bij $A$. So, the 
weight $\psi$ has an obvious extension from $A$ to a faithful, normal 
semi-finite weight on this von Neumann algebra $M$. We will also use 
$\psi$ for this extension.
\snl
It makes sense, and it will be convenient, also to call $\sigma$ the 
modular automorphism group.
\nl
Now we begin with the {\it study of the right invariance} of the weight 
$\psi$ on the C$^*$-algebra $A$ with its comultiplication $\Phi$.
\snl
More precisely, we will prove the following result (we refer to 
definition 2.2 in 
[K-V2] for the notion of invariance that we use here).

\inspr{4.5} Theorem \rm 
For any positive element $x\in A$ 
such that $\psi(x)<\infty$ and any positive linear functional $\omega$ 
in $A^*$, we have that 
$$\psi((\iota\otimes\omega)\Phi(x))=\omega(1)\psi(x)$$
where $\omega(1)=\|\omega\|$. 
\einspr

Recall that slices of 
$\Phi(x)$ belong to $A$ and when we slice with a positive functional, 
we get again something positive. So we can apply $\psi$ and the 
formulation of the invariance above makes sense.
\snl
Before we start with the proof of this theorem, we like to formulate 
two remarks. First, we will show that, in some sense, the classical limit of 
our weight $\psi$ gives the right Haar integral on the classical 
$az+b$-group. This will be a first indication of the fact that we have 
the correct weight. In a second remark, we give some  
argument for the right invariance, not at all precise, but very 
instructive. When we give 
the correct proof afterwards, this will help to understand the steps 
that we take. Moreover, it will show how in fact, this weight can be 
found. We will come back to this in section 5 (where we construct other 
examples) and in section 6.

\inspr{4.6} Remark \rm
Let us consider the classical limit of our system when $q\to 1$.
There is no precise theory for this limit procedure in this context. 
So, neccesarily we will have to be somewhat loose. Let us first 
consider the part with $|a|$ and $|b|$. The weight sends 
$\int f(|b|,t)|a|^{it}\,dt$ to $\int f(r,0) r\,dr$. 
In the limit, where $|a|$ and $|b|$ become commuting elements, a simple 
calculation involving the Fourier transform shows that this results in the 
integral $g\mapsto \frac1{2\pi}\int g(r,s)\frac{r}{s}\,dr\,ds$. 
\snl
On the other hand, let us look at the part with $u$ and $v$. In the 
limit, these unitaries become commuting unitaries with full spectrum. 
The fact that $u^kv^{\ell}$ is sent to $0$ except when $k=\ell=0$ means 
that we will just get a scalar multiple of 
the usual Lebesgue measure on the two-torus $\Bbb 
T^2$. 
\snl
When we take these two parts together (and when we forget about the 
scalars), we see that a function $g$ of two 
complex variables gets the value 
$$\iiiint g(ru,sv) \frac{r}{s} \,dr\,ds\,du
\,dv$$
in terms of polar coordinates and hence
$$\iint g(y,z) \frac{1}{|z|^2} \,dy\,dz$$
in terms of the usual Lebesgue measure on $\Bbb C$. This is precisely 
the right Haar measure on the complex $az+b$-group where the product is 
defined for these variables as $(y,z)(y',z')=(zy'+y',zz')$.
\einspr

In the next remark, we will use right multiplication by elements of the 
C$^*$-algebra $A$ in the GNS-space $\Cal H_\psi$. 
As we will need this notion later 
in this section, let us give a precise definition here.
\snl
Recall that $\Cal N_\psi=\{x\in A\mid \psi(x^*x)<\infty\}$ (by 
definition) and that the map $\eta_\psi$, at first only defined on the 
$^*$-algebra $A_0$, is extended to all of $\Cal N_\psi$ in the process of 
constructing the weight.

\inspr{4.7} Definition \rm For elements $x\in A$, say analytic with 
respect to the modular automorphism group $\sigma$, we define the bounded 
operator $\pi_\psi'(x)$ by $\pi_\psi'(x)\eta_\psi(y)=\eta_\psi(yx)$ 
whenever $y\in \Cal N_\psi$.
\einspr

In fact, by a standard argument, we have 
$\pi'_\psi(x)=J\pi_\psi(\sigma_{\frac{i}{2}}(x))^*J$
and this formula holds for all elements in the domain of 
$\sigma_{\frac{i}{2}}$.
Of course, $\pi_\psi'$ is a anti-representation and not a 
$^*$-representation (and it is unbounded). 

\inspr{4.8} Remark \rm
Let us suppose for a moment 
that $\psi$ is the correct right invariant weight. Then we can 
consider the right regular representation. It is a unitary $\widetilde 
W$ on $\Cal H_\psi\otimes\Cal H_\psi$. If we use the Sweedler notation
$\Phi(x)=\sum x_{(1)}\otimes x_{(2)}$, then we can write formally 
$$\widetilde W(\eta_\psi(x)\otimes\xi)=\sum \eta_\psi(x_{(1)})\otimes 
\pi_\psi(x_{(2)})\xi$$
whenever $x\in A_0$ and $\xi\in\Cal H_\psi$.
\snl
Now, observe that $\Phi(x)=W(x\otimes 1)W^*$. And recall that in this very 
special case, $W\in M\otimes M$ where $M$ is the associated von Neumann 
algebra (see section 3). This implies that, still formally, we can also apply 
$\pi_\psi$, as well as $\pi'_\psi$ on the first leg of $W$. Then we can 
write $\widetilde W=W_1W_2$ where
$$\align W_1&=(\pi_\psi\otimes \pi_\psi)W\\
         W_2&=(\pi'_\psi\otimes\pi_\psi)W^*.
\endalign$$
The first factor $W_1$ presents no problem. It is a unitary as $\pi_\psi$ is 
a $^*$- representation. The second factor $W_2$ is more tricky. As  
$W^*=(\hat S\otimes\iota)W$ (cf.\ 3.11), we have that $W_2$ is
of the form  $(\pi'_\psi\circ\hat S\otimes \pi_\psi)W$.
Now, we will show (see proposition 4.10 below), 
that $\pi'_\psi\circ\hat S$ is in fact also a 
$^*$-representation (and bounded), 
now from $M$ to $\pi_\psi(M)'$. Therefore, the 
second factor is also a unitary. Finally observe that the unitarity of 
$\widetilde W$ is essentially equivalent with the right invariance of 
$\psi$ (see the proof of the theorem below).
\einspr

So we see that the key is that the 
action of the modular automorphism group compensates the result of the scaling 
group. These automorphism groups are the obstructions for $\pi'_\psi$ 
and $\hat S$ to be $^*$-maps.
\snl
By the structure of the von Neumann algebra $M$, the weight is 
determined by the modular automorphism group. Using all these ideas, it 
is possible to construct this weight. This is not the way how we found 
it, but we could have used this method. We will say something more 
about this in sectin 5, where we give other examples. 
We will also give some more comments in section 6 
where we make futher remarks and where we refer to forthcoming papers 
about this procedure.
\snl
As we see from the preceding remark, in order to prove the invariance 
of $\psi$, we will essentially need to show that $\pi'_\psi\circ\hat S$ 
is a $^*$-representation. This will be done in proposition 4.10. We 
first need the following lemma. It is a key result: it is 
really this property that eventually gives the right invariance of 
$\psi$. 

\inspr{4.9} Lemma \rm
Let $(a,b)$ be an admissible pair of normal operators as before. Let 
$(\hat a,\hat b)$ be the pair associated to $(a,b)$ as in 3.9, given by
$\hat a=b^{-1}$ and $\hat b=ab^{-1}$. Consider the scaling group 
$\hat\tau$ as defined in proposition 3.12, but for the pair $(\hat 
a,\hat b)$. So, $\hat\tau_t(\hat b)=e^{\frac{2\pi t}{n}}\hat b$
and $\hat\tau_t(\hat a)=\hat a$ for all $t$. Then $\hat\tau$ 
coincides with the modular automorphism $\sigma$ as defined in the 
formulation of the theorem 4.4.

\snl\bf Proof: \rm
We have 
$$\align
\sigma_t(\hat a)&=\sigma_t(b^{-1})=b^{-1}=\hat a = \hat\tau_t(\hat a) \\ 
\sigma_t(\hat b)&=\sigma_t(ab^{-1})=\sigma_t(a)\sigma_t(b^{-1})\\
&=e^{\frac{2\pi t}{n}} ab^{-1}=e^{\frac{2\pi t}{n}}\hat b=
\hat\tau_t(\hat b).
\endalign$$
\einspr

We see that the proof is essentially trivial. But it is important to 
have the precise statement. It might be confusing because we are working 
with two pairs $(a,b)$ and $(\hat a,\hat b)$, closely related but 
playing a different role.
\snl
An immediate consequence of this property is what we need:

\iinspr{4.10} Proposition \rm Let $(a,b)$ and $(\hat a,\hat b)$ be as 
before. Consider the antipode $\hat S$ as defined for the pair 
$(\hat a,\hat b)$ 
and $\pi'_\psi$ as defined in the preceding for $(a,b)$. Then 
$\pi'_\psi\circ \hat S$ extends to an injective normal $^*$-homomorphism 
$\gamma$ of the von Neumann algebra $M$ to $\pi_\psi(M)'$.

\snl\bf Proof: \rm
We could prove this by verifying it on the generators. This would give 
necessary formulas for later results in the paper.
\snl 
However, we want to use the previous result. Indeed,
when $x$ is in the domain of $\hat\tau_{-\frac{i}{2}}$, then
$$\align 
(\pi'_\psi\circ\hat S)(x)
&=J|T|\pi_\psi(\hat S(x))^*|T|^{-1}J \\
&= J\sigma_{-\frac{i}{2}}(\pi_\psi(\hat R\hat 
\tau_{-\frac{i}{2}}(x))^*)J \\
&= J\sigma_{-\frac{i}{2}}(\pi_\psi(\hat\tau_{\frac{i}{2}}(\hat R(x))^*)J \\
&= J\pi_\psi(\hat R(x^*))J.
\endalign$$
So we can define $\gamma:M\to \pi_\psi(M)'$ by 
$\gamma(x)=J\pi_\psi(\hat R(x^*))J$.
\einspr

Later in this section, we will give an explicit formula for the right 
regular representation $\widetilde W$ which, as we just saw, should be
$W_1W_2=(\pi_\psi\ot\pi_\psi)W(\gamma \ot \pi_\psi)W$.
\snl
Now we are ready to prove the invariance.

\inspr{} Proof \rm of theorem 4.5: 
Take any $x\in A^+$ satisfying $\psi(x)<\infty$ and an element $\omega\in 
A^*$ such that $\omega\geq0$. We must consider 
$(\iota\otimes\omega)\Phi(x)$. We use that $\Phi(x)=W(x\otimes 1)W^*$ 
where $W$ is the multiplicative unitary as in proposition 3.10 of the 
previous section. Without loss of generality, we can assume that the 
second leg of $W$ is faithfully represented on a Hilbert space $\Cal K$ 
where $\omega$ is a vector state. So we can assume that 
there is a vector $\xi\in\Cal K$ 
such that $\omega=\omega_{\xi,\xi}$. We use $\omega_{\xi_1,\xi_2}$ to 
denote the linear functional on $\Cal B(\Cal K)$ defined by 
$\omega_{\xi_1,\xi_2}(x)=\langle x \xi_1, \xi_2\rangle$.
\snl
Now let $(\xi_i)$ be an orthonormal basis in $\Cal K$. Then we can 
write
$$\align
(\iota\ot\omega)\Phi(x)
&=(\iota\ot\omega_{\xi,\xi})(W(x\ot 1)W^*)\\
&=\sum_i ((\iota\ot \omega_{\xi_i,\xi})W)x
((\iota\ot \omega_{\xi,\xi_i})W^*).
\endalign$$
where the sum is convergent in the $\sigma$-weak topology on the von 
Neumann algebra $M$.
\snl
We know that $(\iota\ot\omega_{\xi,\xi_i})W$
belongs to the domain of $\hat S$ and we have 
$\hat S((\iota\ot \omega_{\xi,\xi_i})W)=(\iota\ot\omega_{\xi,\xi_i})W^*$ 
(see lemma 3.11). So, this element also belongs to the domain of
$\hat\tau_{-\frac{i}{2}}$. As $\hat\tau$ coincides with $\sigma$, we 
have also that this element belongs to the domain of 
$\sigma_{-\frac{i}{2}}$. Then it follows that 
$(\iota\ot\omega_{\xi,\xi_i})W^*$ is in the domain of 
$\sigma_{\frac{i}{2}}$. Hence, it is right 
bounded and we get 
$x^{\frac12}(\iota\ot \omega_{\xi,\xi_i})W^* \in \Cal N_\psi$ and
$$\eta_\psi(x^{\frac12}(\iota\ot \omega_{\xi,\xi_i})W^*)
=\gamma((\iota\ot\omega_{\xi,\xi_i})W)\eta_\psi (x^{\frac12})$$
where $\gamma$ is the 
$^*$-representation that we obtained in proposition 4.10.
If we write $W_2=(\gamma\ot\iota)W$, 
then we know that $W_2$ is 
unitary. If we now use that $\psi$ is a normal weight, we get
$$\align \psi((\iota\ot\omega)\Phi(x))&=\sum \|((\iota\ot
\omega_{\xi,\xi_i})W_2) \eta_\psi(x^\frac12)\|^2\\
&=\|W_2(\eta_\psi(x^\frac12)\ot\xi)\|^2\\
&=\psi(x)\langle \xi,\xi\rangle = \omega(1)\psi(x).
\endalign$$
This proves the result.
\einspr
 
In general, to prove the invariance of a weight can be rather hard. 
Here it turns out to be relatively simple. One reason is that, 
according to the general theory, we only have to consider positive elements 
$x\in A$ such that $\psi(x)<\infty$. On the other hand, we are using a 
special and useful technique here. One might think that this will only 
work in very special cases (as the example here and the ones that we 
treat in section 5). However, as we will explain in section 6, there 
are good reasons to believe that it will also work in many more cases. 
\snl
Now we have essentially proven the following result.

\iinspr{4.11} Theorem \rm 
The pair $(A,\Phi)$ is a locally compact quantum group.
\einspr

Indeed, from the existence of the unitary antipode $R$, which is a 
$^*$-anti-isomorphism of $A$ that flips the comultiplication, we also 
find the existence of a suitable left invariant weight $\varphi$ 
defined by $\psi\circ R$. The density conditions needed to have a 
locally compact quantum group have already been discussed in section 3). 
So all the axioms are fulfilled and the theorem is proven.
\snl
The left invariant weight $\varphi$ can be constructed by using the 
formula $\psi\circ R$. This however, is not so simple. It can also be 
constructed by other methods. We plan to include the explicit form of 
the left invariant weight in a later version of 
this paper.
\snl

Now, it is easy to see that the Haar weights are not invariant with respect 
to the scaling group for this locally compact quantum group. As we have 
explained already, this is an important new feature. We obtain 
the concrete scaling factor in the following proposition.

\iinspr{4.12} Proposition \rm
The right invariant weight $\psi$ is relatively invariant with respect 
to the scaling group. More precisely, we have
$$\psi(\tau_t(x))=e^{-\frac{4\pi t}{n}}\psi(x)$$
for all $x\in A^+$.

\snl\bf Proof: \rm
Recall that
$$\psi(x)=\int f_{0,0}(r,0)r\,dr$$
when 
$$x=\sum_{k,\ell}\left(\int 
f_{k,\ell}(|b|,t)|a|^{it}\,dt\right)v^ku^\ell$$
as in 3.5. Now, by 3.12, the scaling group $\tau$ has the property that 
$\tau_s(b)=e^{\frac{2\pi s}{n}}b$ and $\tau_s(a)=a$. Therefore
$$\tau_s(x)=\sum_{k,\ell}\left(\int 
f_{k,\ell}(e^{\frac{2\pi s}{n}}|b|,t)|a|^{it}\,dt\right)v^ku^\ell$$
and so
$$\align
\psi(\tau_s(x))&=\int f_{0,0}(e^{\frac{2\pi s}{n}}r,0)r\,dr\\
&=e^{-\frac{4\pi s}{n}}\int f_{0,0}(r,0)r\,dr\\
&=e^{-\frac{4\pi s}{n}}\psi(x).
\endalign$$
Now, it will follow easily that this gives the appropriate scaling on the 
Hilbert algebra level and hence also the full weight on the von Neumann 
algebra and on the C$^*$-algebra will have the same scaling property.
\einspr

The equality of $\psi\circ\tau_s$ and $e^{\frac{-4\pi s}{n}}\psi$ on the 
dense $^*$-subalgebra $A_0$ yields the overall equality of these two 
weights. This follows from the construction as we have argued at the 
end of the proof. Of course, it can also be shown using other standard 
techniques. Observe that the modular automorphism group commutes with 
the scaling group. This is easily verified here because we have 
explicit formulas for these automorphisms in terms of the generators. 
In fact, this also is a general result (see proposition 6.8 in [K-V2]).
\nl 
Finally, we will give an explicit formula for the right regular 
representation. It is a multiplicative unitary, similar to the original 
one, but different. We will show how the one is related with the other.
\snl
First, we need the explicit formulas for $\pi'_\psi$ on the generators. 
We give them in the following lemma.

\iinspr{4.13} Lemma \rm
We have
$$\align
\pi'_\psi(a)  &=q(1 \ot a_1 \ot 1 \ot s) \\
\pi'_\psi(a^*)&=q(1 \ot a_1 \ot 1 \ot s^*) \\
\pi'_\psi(b)  &=b_0 \ot c_1 \ot s \ot m \\
\pi'_\psi(b^*)&=b_0 \ot c_1 \ot s^* \ot m^* .
\endalign$$
\einspr

Observe again that $\pi'_\psi$ is a anti-representation but not a 
$^*$-anti-representation, the obstruction being related with the 
modular automorphisms. This is reflected by the fact that 
$\pi'_\psi(a^*)$ is not the adjoint of $\pi'_\psi(a)$ but rather 
 $q^2\pi'_\psi(a)^*$. This is different with the other generator where 
we do have $\pi'_\psi(b^*)=\pi'_\psi(b)^*$.
\snl
Recall that 
$\sigma_t(x)=|T|^{2it}x|T|^{-2it}$ and that
$$\align
\sigma_t(\pi_\psi(a))&=e^{\frac{2\pi t}{n}}\pi_\psi(a)\\
\sigma_t(\pi_\psi(b))&=\pi_\psi(b).
\endalign$$
In particular, $\sigma_t$ only scales $\pi_\psi(|a|)$ and leaves the 
other generators $u$, $|b|$ and $v$ invariant. This also clarifies the 
obstruction factor $q$ in the formulas for $\pi'_\psi(a)$ and 
$\pi'_\psi(a^*)$.
\snl
It is also instructive to verify the formulas
$$J|T|\pi_\psi(x)=\pi'_\psi(x^*)J|T|,$$
at least formally, for the different generators $a$, $a^*$, $b$ and 
$b^*$ (or equivalently on $u$, $|a|$, $v$ and $|b|$).

\iinspr{4.14} Proposition \rm
The regular representation $\widetilde W$ (see remark 4.8) has the 
following expression
$$\widetilde W=F(\hat b\ot b)\chi(\hat a\ot 1,1\ot a)$$
where $F$ and $\chi$ are the special functions as desribed in the 
previous section (before the definition 3.7) 
and where now $(a,b)$ and $(\hat a,\hat b)$ 
are the admissible pairs of normal operators on the space $\Cal 
H_\psi=L^2(\Bbb R^+)\ot L^2(\Bbb R)\ot \Bbb C^{2n}\ot \Bbb C^{2n}$
given by
$$\align
a &=a_0\ot a_1 \ot m \ot s \\
b &=b_0\ot 1 \ot s \ot 1 \\
\hat a &=1\ot c_1 \ot 1 \ot m \\
\hat b &=\,\cdot\,\ot a_1 \ot \,\cdot\, \ot s 
\endalign$$
with the part $\,\cdot\,\ot \,\cdot\,$ on $L^2(\Bbb R^+)\ot \Bbb C^{2n}$
given by (the closure of)
$$(a_0b_0^{-1}\ot ms^*)-q^{-1}(b_0^{-1}\ot s^*).$$

\snl\bf Proof: \rm
From the previous observations (in particular remark 4.8 and the remark 
following proposition 4.10), the definition of $W$ in 3.7 and the 
formulas in proposition 4.2, 
we know that $\widetilde W=W_1W_2$ where
$$\align 
W_1 &= F(\hat b_1\ot b)\chi(\hat a_1\ot 1,1\ot a)\\
W_2 &= F(\hat b_2\ot b)\chi(\hat a_2\ot 1,1\ot a)
\endalign$$
with $a$ and $b$ as in the formulation of the proposition and $\hat 
a_1=b^{-1}$ and $\hat b_1=ab^{-1}$ and $\hat a_2=\gamma(\hat a_1)$ and 
$\hat b_2=\gamma (\hat b_1)$. If moreover, we use the definition of 
$\gamma$ as $\pi'_\psi\circ \hat S$ and the formulas in the previous 
lemma, we get
$$\align
\hat a_2 &=\pi'_\psi(\hat S(\hat a))= \pi'_\psi(\hat a^{-1}))
    =\pi'_\psi(b)=b_0\ot c_1 \ot s \ot m \\
\hat b_2 &=\pi'_\psi(\hat S(\hat b))=\pi'_\psi(-\hat a^{-1}\hat 
b)=-\pi'_\psi(bab^{-1})\\
         &=-q(1\ot c_1^{-1}a_1c_1 \ot 1 \ot m^*sm)
	 =-q^{-1}(1\ot a_1 \ot 1 \ot s).
\endalign$$
Now, the bicharacter $\chi$ has the property that 
$$\chi(\gamma,a)b\chi(\gamma,a)^*=\gamma b$$
(cf.\ formula 2.2 in [W5]) because $(a,b)$ is an admissible pair. 
Observe that $\hat a_1$ and $\hat b_2$ are commuting operators. It 
follows that
$$\align \chi(\hat a_1 \ot 1, 1\ot a)&(\hat b_2\ot b) 
\chi(\hat a_1 \ot 1, 1\ot a)^*\\	 	   
&=(\hat b_2\ot 1)\chi(\hat a_1 \ot 1, 1\ot a)(1\ot b)
\chi(\hat a_1 \ot 1, 1\ot a)^*\\
&=(\hat b_2\ot 1)(\hat a_1\ot b)=\hat b_2\hat a_1\ot b.
\endalign$$
Therefore
$$\widetilde W=W_1W_2=F(\hat b_1\ot b)F(\hat b_2\hat a_1\ot b)
\chi(\hat a_1 \ot 1, 1\ot a)\chi(\hat a_2 \ot 1, 1\ot a).$$
Now, we can use the exponential properties. Observe that
$$\align
\hat b_1&=ab^{-1}=a_0b_0^{-1}\ot a_1 \ot ms^* \ot s \\
\hat b_2\hat a_1&=-q^{-1}b_0^{-1}\ot a_1 \ot s^* \ot s.
 \endalign$$
This means that the exponential formula (see [W5]) can be used. And 
$$F(\hat b_1\ot b)F(\hat b_2\hat a_1\ot b)=F(\hat b\ot b)$$
with $\hat b$ as in the formulation of the proposition.
For the second part, we have
$$\chi(\hat a_1 \ot 1, 1\ot a)\chi(\hat a_2 \ot 1, 1\ot a)=
\chi(\hat a_1\hat a_2 \ot 1, 1\ot a)$$
and indeed $\hat a_1\hat a_2= 1\ot c_1\ot 1 \ot m$ and this is how we 
defined $\hat a$ here.
\einspr

In a later version of this paper, we will make 
a comparison of the new multiplicative 
unitary $\widetilde W$ and the originial one $W$ as defined in 3.7.
\nl
\nl



\bf 5. Other examples \rm
\nl

In this section, we will treat some other examples. 
The first one is very similar to the one that we have already given in 
full detail. It is also a quantization of the $az+b$-group, but now 
with a real deformation parameter $q$. This example is briefly 
considered by Woronowicz in an appendix in [W5]. For this example, the 
Haar weights are invariant with respect to the scaling group. There is 
also the quantization of the $ax+b$-group as studied by Woronowicz in 
[W-Z]. Here again the Haar measures are  not invariant but only 
relatively invariant for the scaling group. As we mentioned already in 
the introduction, it was this example that we discovered first having 
this non-invariance property.
\snl
We will not treat these examples in full detail but only mention those 
aspects that are different from the complex deformation of the 
$az+b$-group that we studied in this paper. At the end of the section, 
we will briefly discuss the quantum $E(2)$ and its dual. The Haar 
measures on these quantum groups were already obtained before in [B1] 
and [B2].
\nl 
So we start with the quantization of the $az+b$-group with a real 
deformation parameter $q$. Now $q$ is supposed to satisfy $0<q<1$. 
The underlying Hopf 
$^*$-algebra is the one of proposition 2.2 (with $\lambda=q^2$ as 
before). But because $q$ is real, the commutation rules translate 
into different commutation rules for the elements in the polar 
decomposition. Moreover, the spectral conditions become of a different 
nature. Whereas in the complex case, it was possible to impose them in 
an algebraic way (by requiring $a^n$ and $b^n$ to be self-adjoint like 
in proposition 2.3), this is no longer the case here.
\snl
Here is how definition 3.1 has to be adapted (see appendix A of [W5]).

\inspr{5.1} Definition \rm
Let $(a,b)$ be a pair of normal operators on a 
Hilbert space $\Cal H$. Let $a=u|a|$ and $b=v|b|$ be the polar 
decompositions of $a$ and $b$. Assume that $|a|$ and $|b|$ are 
non-singular so that $u$ and $v$ are unitary. Furthermore, assume that 
\snl
i) the spectra of $|a|$ and $|b|$ are contained in the set 
$\{q^n\mid n\in \Bbb Z \}\cup \{0\}$,

ii) $|a|v=qv|a|$ and $|b|u=q^{-1}u|b|$,

iii) $uv=vu$,

iv) $|a|$ and $|b|$ strongly commute.
\snl
Then we call $(a,b)$ and {\it admissible pair of normal operators}.
\einspr
 
If we compare this with definition 3.1, we see that now the spectra of 
$a$ and $b$ are restricted to the closure $\overline\Gamma$ of the 
group $\Gamma$ defined as
$$\Gamma = \{zq^n \mid n\in \Bbb Z, z\in \Bbb C \text{ and } |z|=1 \}$$
Here $\Gamma$ is the group $\Bbb T \times \Bbb Z$ and again it is 
self-dual. However, the self-duality is of a different nature as 
in the complex case. There, the group was a direct product of two 
self-dual factors while here, the duality interchanges the factors. 
There are reasons to believe that it is this difference between the two 
cases that is responsible for the fact that the scaling group leaves 
the Haar measures invariant in the second case and not in the first 
case.
\snl
Also here, it is not so difficult to construct elementary representations. 
As before, they will turn out to be the building blocks for the other 
representations we will have later.

\inspr{5.2} Proposition \rm
 Consider the Hilbert space $\ell^2(\Bbb Z)$ with an 
orthonormal basis $(e_k)_{k\in \Bbb Z}$. Define operators $m$ and $s$ 
by
$$\align me_k &= q^k e_k \\ se_k & =e_{k+1}.
\endalign$$ 
Then, the operators on $\ell^2(\Bbb Z)\ot \ell^2(\Bbb Z)$ given by
$$\align a &= m\ot s^* \\ b &= s\ot m
 \endalign$$
satisfy the conditions in definition 5.1.
\einspr

The basic commutation rule $ms=qsm$ gives all the other ones. Observe 
that the polar decompositions are as follows:
$$\matrix |a|=m\ot 1 &\qquad\qquad\qquad & u=1\ot s^* \\
          |b|=1\ot m &\qquad\qquad\qquad & v=s\ot 1. 
\endmatrix$$ 
\snl
Again it follows from the general theory (cf.\ the appendix, example 
A.4.vi) that any 
admissible pair of normal operators $(a,b)$ 
(that is satisfying the conditions of definition 5.1), 
is obtained from this irreducible pair by tensoring with 
$1$ on some Hilbert space.
\nl
Now, we come to the {\it C$^*$-algebra}. Given the commutation rules in 5.1, 
it is quite clear that we must have something like the following:

\inspr{5.3} Proposition  \rm Consider an admissible pair of operators 
on $\Cal H$ 
as in 5.1. Let $A_0$ be the space of bounded linear operators of the 
form 
$$\sum_{k,\ell} f_{k,\ell}(|a|,|b|)v^ku^\ell$$
where $f_{k,\ell}$ are functions with finite support of two variables 
in $\{q^n\mid n\in \Bbb Z\}$ and such that only finitely many 
$f_{k,\ell}$ are non-zero. 
Then $A_0$ is a $^*$-algebra, acting non-degenerately on $\Cal H$.
\einspr

The proof of this result is straightforward. 
When $x$ is the operator above, then
$$x^*=\sum_{k,\ell} \overline{f_{k,\ell}}(q^k|a|,q^{-\ell}|b|)v^{-k}u^{-\ell}$$
and
$$x^*x=\sum_{k,\ell,k',\ell'} \overline{ f_{k,\ell}}(q^k|a|,q^{-\ell}|b|)
f_{k',\ell'}(q^k|a|,q^{-\ell}|b|) v^{k'-k}u^{\ell'-\ell}.$$
We will need these formulas later.
\snl
The C$^*$-algebra that we need is somewhat bigger than the norm closure 
of this $^*$-algebra. We must take the C$^*$-algebra 'generated' by $a$, 
$a^{-1}$ 
and $b$ (in the appropriate sense). Therefore, we must 
allow functions $f_{k,\ell}$ that have 
a non-zero limit when the second variable tends to $0$.  
We will use $A$ for this bigger C$^*$-algebra. As  before, it 
contains (an isomorphic image of) the crossed product of
$C_0(\Bbb Z\times\Bbb Z)$ by the action $\alpha$ of $\Bbb Z\times\Bbb Z$ given 
by $(\alpha_{k,\ell}g)(r,s)=g(r-k,s+\ell)$.  
\snl
For the von Neumann algebra $M$, we can just take the weak closure of 
the $^*$-algebra $A_0$ and we don't have to worry about these 
restrictions. We have mentioned before already that the von Neumann 
framework is easier and that, only from a theoretical point of view, it 
can make sense to consider the C$^*$-framework in concrete examples.
\snl
Finally, observe that both the C$^*$-algebra 
$A$ and the von Neumann algebra $M$ do not depend on the choice of the 
pair $(a,b)$ (see again the appendix).
\nl
The {\it comultiplication} $\Phi$ on this C$^*$-algebra $A$ and on the von 
Neumann algebra $M$ is again described by 
a multiplicative unitary $W$ of the form 
$$F(\hat b\ot b)\chi(\hat a\ot 1,1\ot a)$$
where $\chi$ is the appropriate bicharacter expressing the self-duality 
of the underlying group $\Gamma$ and where is $F$ is the appropriate 
version of the quantum exponential function. Now the convention, used 
by Woronowicz here, is 
$\hat a=b^{-1}$ and  $\hat b=b^{-1}a$. This is different from the 
complex case (see lemma 3.9 and proposition 3.10). 
There is no need to go into details here for what we need. So, we refer to 
[W5].
\snl
However, we need the formula for the scaling group. It is given in the 
following proposition.

\inspr{5.4} Proposition \rm
There exists a strongly continuous one-parameter group of 
$^*$-automor\-phisms $(\tau_t)$ of $A$ defined by $\tau_t(a)=a$ and 
$\tau_t(b)=q^{-2it}b$. And there is an involutive 
$^*$-anti-automorphism $R$ of $A$ given by $R(a)=a^{-1}$ and 
$R(b)=-qa^{-1}b$. Together, they give the polar decomposition of the 
antipode (as in proposition 3.12).
\einspr

In this case, it is easier to show that $R$ is well-defined. The polar 
decomposition of $qa^{-1}b$ is $u^*v|a|^{-1}|b|$ and if we consider the 
representation given in 5.2, we see that 
$$\align 
u^*v &=s\ot s \\ 
|a|^{-1}|b| &=m^{-1}\ot m
\endalign$$
and it is standard to construct a unitary $U$ that commutes with $u$ and 
$|a|$ and that transforms $s \ot 1$ into $s \ot s$ 
and $1\ot m$ into $m^{-1}\ot 
m$. In fact, it is given by $U(e_k\ot e_\ell) =e_k\otimes e_{k+\ell}$.
This will take care of the non-trivial part of $R$.
\nl
The next step is the construction of the {\it right Haar measure}.

\inspr{5.5} Proposition \rm
There exists a lower semi-continuous densely defined, faithful KMS-weight 
$\psi$ on $A$ given by 
$$\psi(x)=\sum_{i,j}f_{0,0}(q^i,q^j)q^{2j}$$
when $x$ has the form 
$$x=\sum_{k,\ell} f_{k,\ell}(|a|,|b|)v^ku^\ell$$
as in 5.3.
\einspr

Remark that indeed, the fact that we also consider functions $f_{0,0}$ 
with a possible non-zero limit when the second variable goes to $0$ 
presents no problem for the convergence in this sum as $0<q<1$.
\snl
The proof is standard (and similar as in section 4). From
$$\align\psi(x^*x)&=\sum_{k,\ell,i,j} 
|f_{k,\ell}(q^kq^i,q^{-\ell}q^j)|^2q^{2j}\\
&=\sum_{k,\ell,i,j} 
|f_{k,\ell}(q^i,q^j)|^2q^{2j}q^{2\ell}
\endalign$$
we see that we can identify the GNS-space with the tensor product of 
four copies of $\ell^2(\Bbb Z)$ and then the canonical map 
$\eta_\psi$ is given by 
$$\eta_\psi(x)=\xi=\sum\xi_{k,\ell}\ot e_k\ot e_\ell$$
where
$$\xi_{k,\ell}(i,j)=q^{\ell+j}f_{k,\ell}(q^i,q^j).$$
It is also straightforward to calculate the different operators 
involved. We have, using the notations of 5.2 (and the notations $\pi$ 
and $\pi'$ from section 4)
$$\align 
\pi(a)&=m\ot s^*\ot 1\ot s \\
\pi(b)&=s\ot m  \ot s \ot1\endalign$$
and
$$\align
\pi'(a)&=q(m\ot 1\ot m^{-1}\ot s) \\
\pi'(a^*)&=q^{-1}(m\ot 1\ot m^{-1}\ot \ot s^*) \\
\pi'(b)&= 1\ot m\ot s\ot m \\
\pi'(b^*)&=1\ot m\ot s^*\ot m.
\endalign$$ 
Observe that the 'obstruction' for $\pi'$ to be a 
$^*$-anti-automorphism lies in $\pi'(u)=q(1\ot 1\ot 1\ot s)$ which is 
not a unitary but a scalar multiple of a unitary. This is in fact 
already an indication that in this example the Haar weights will be 
invariant w.r.t.\ the scaling group (as we will see later).
\snl
And as before, this is related with the action of the modular automorphisms. 
One finds
$$T(e_{i}\ot e_{j}\ot e_{k}\ot e_{\ell})=
q^{-\ell}(e_{i-k}\ot e_{j+\ell}\ot e_{-k}\ot e_{-\ell})$$
and so
$$|T|=1\ot 1\ot 1\ot m^{-1}.$$
This operator commutes with $|a|$, $|b|$ and $v$, but not with $u$. In 
fact, when as before we use $\sigma_t$ also for the automorphism on $A$ 
induced by the modular automorphism $\sigma_t=|T|^{2it}\,\cdot\,|T|^{-2it}$
on $\pi(A)$, we find
that $\sigma_t(a)=q^{-2it}a$ and $\sigma_t(b)=b$.
\snl
As we know from the discussion in section 4, the basic 
result to prove the right invariance is that the modular automorphism 
group $\sigma$ coincides with the scaling group $\hat \tau$. This will be done 
for this example in proposition 5.6.
\snl
Before we do this however, let us 'verify' the classical limit and see 
what happens when $q$ approaches $1$. If we write
$$\sum f(q^i,q^j)q^{2j}=\sum f(q^i,q^j)
\frac{1}{q^i}\left(\frac{q^{i+1}-q^i}{(q-1)}\right)
q^j\left(\frac{q^{j+1}-q^j}{q-1}\right) 
$$
we see that in the limit $q\to 1$, the expression 
$(q-1)^2\sum f(q^i,q^j)q^{2j}$
will precisely transform into
$$\iint f(r,s)\frac{1}{r}dr\, sds.$$
Together with the unitary part, we will arrive at the functional $\psi$ 
given by
$$\psi(f)=\iiiint f(re^{i\varphi},se^{i\theta}) \frac{1}{r}dr\, 
sds\,d\varphi d\theta$$
and this is precisely the right invariant integral on the classical 
$az+b$-group as we saw in section 4. Observe that the variables $a$ and 
$b$ are interchanged.
\snl
This is one reason why we can expect to have the correct right 
invariant functional on this quantization of the $az+b$-group. Of 
course, the real reason is the following result.

\inspr{5.6} Proposition \rm
Let $(a,b)$ be an admissible pair of normal operators (as in 5.1). Let 
$(\hat a,\hat b)$ be the associated pair given by $\hat a=b^{-1}$ and 
$\hat b=b^{-1}a$. Then, the scaling group $\hat \tau$ as defined in 
proposition 5.4, but now for the pair   $(\hat a,\hat b)$, coincides 
with the modular automorphism group $\sigma$.

\snl\bf Proof: \rm
$$\align
\sigma_t(\hat a)&= \sigma_t(b^{-1})=b^{-1}= \hat a=\hat\tau_t(\hat a) \\
\sigma_t(\hat b)&=\sigma_t( b^{-1}a)=q^{-2it}b^{-1}a=q^{-2it}\hat 
b=\hat\tau_t(\hat b).
\endalign$$ 
\einspr

Now, the argument continues as in section 4 for the complex case.
\snl
Here however, as we mentioned already, the Haar weight is invariant 
with respect to the scaling group. Because $\tau_t$ leaves $a$ 
invariant, it will also leave $u$ and $|a|$ invariant. And as 
$\tau_t(b)=q^{-2it}b$ we have $\tau_t(v)=q^{-2it}v$ and 
$\tau_t(|b|)=|b|$. It follows easily that $\psi(\tau_t(x))=\psi(x)$ for 
all $t\in\Bbb R$ when $x$ is in $A_0$ as in 5.5.
\snl
Again, using the previous formulas, the commutation rules and the 
exponential properties of the special functions involved, one can 
calculate the regular representation. It has the form $\widetilde 
W=F(\hat b\ot b)\chi(\hat a\ot 1,1\ot a)$ where now 
$$\align \hat a &= s^*\ot 1\ot 1\ot m \\
         \hat b &=s^*m \ot \,\cdot\, \ot \,\cdot\, \ot s
\endalign$$
and where the part $\,\cdot\,\ot \,\cdot\,$ acts on 
$\ell^2(\Bbb Z)\ot\ell^2(\Bbb Z)$ as
$$ -m^{-1}\ot m^{-1}s^*+ m^{-1}s^* \ot s^*.$$
(where of course, the closure of the sum has to be taken).
\nl\nl     
The {\it second example} that we will treat in this section, is the 
quantized $ax+b$-group.
\snl
Here, the starting point is the Hopf $^*$-algebra obtained from the Hopf 
algebra of proposition 2.1 and requiring that $a$ and $b$ are 
self-adjoint. This restricts the deformation parameter $\lambda$ to 
$|\lambda|=1$ as we explained already after this proposition. 
\snl
It is possible to associate a natural C$^*$-algebra, but it seems to be 
impossible to lift the comultiplication to the C$^*$-level (see e.g.\ 
[W-Z]). Woronowicz and Zakrzewski were 
able to solve this problem by adding an extra 
generator. Unfortunately, the comultiplication applied to this 
generator has no simple expression and this extension cannot be 
formulated on the Hopf $^*$-algebra level.
\snl
We refer to [W-Z] for more information about this difficulty and how to 
overcome it. What we will do here is start with the operator 
realization of the generators, just as we did for the other examples in 
3.1 and 5.1. 
\snl
We fix a real number $\theta$ in $\R$ and assume that $0<\theta<\pi$.

\inspr{5.7} Definition \rm
Let $a$ and $b$ be self-adjoint operators on a Hilbert space $\Cal H$. 
Assume that $a$ is non-singular and positive. Let $b=v|b|$ be the polar 
decomposition of $b$. Assume that also $b$ is non-singular so that $v$ 
is a self-adjoint element satisfying $v^2=1$ that commutes with $|b|$. Finally
let $w$ be another self-adjoint element satisfying $w^2=1$. We will 
call $(a,b,w)$ and {\it admissible triple} if also the following 
conditions are satisfied:
\snl
i) $a^{it}|b|a^{-it}=e^{t\theta}|b|$ for all $t\in \R$,

ii) $vw=-wv$,

iii) $v$ and $w$ commute with $a$ and $|b|$.

\einspr 

This approach is a little different from the one of Woronowicz in [W-Z]. 
It is closer in spirit to the other two cases we had already (observe 
e.g.\ the similarity of definition 5.7 with definition 3.1). This point 
of view is also more convenient for our purpose. In this case, the 
underlying group is $\R\times\Z$, but the role played by this group is 
not of the same type as in the two other examples.
\snl
Observe that we are using a different notation than in [W-Z]. There 
$\theta$ is denoted by $\hbar$ while $\beta$ is used instead of $w$.
\snl
It is not difficult to construct the obvious irreducible representation 
here.

\inspr{5.8} Proposition \rm 
Consider the Hilbert space $L^2(\R^+)$ where $\R^+$ is 
considered with the usual Lebesgue measure. Define self-adjoint 
positive non-singular operators $a_0$ and $b_0$ on $L^2(\R^+)$ by
$$\align (a_0^{is}f)(u)&=e^{\frac12 s\theta}f(e^{s\theta} u) \\
         (b_0 f)(u)&=uf(u)
\endalign$$
where $u\in \R^+$ and $s\in\R$. Also consider $\C^2$ with a basis 
$(e_0,e_1)$, indexed over the group $\Z_2$, and define operators $m$ 
and $s$ given by 
$$\align me_k &=(-1)^ke_k\\ se_k &= e_{k+1}. \endalign$$
If we let $a=a_0 \ot 1$, $b=b_0 \ot m$ and $w=1\ot s$, we get operators 
satisfying the properties of definition 5.7.
\einspr

The C$^*$-algebra taken here is the one given in the following 
proposition.

\inspr{5.9} Proposition \rm Consider a triple $(a,b,w)$ as in 
definition 5.7 and let $b=v|b|$ be the polar decomposition of $b$ as 
before. Let $A_0$ be the space of bounded linear operators of the form
$$x=\sum_{k,\ell=0}^{1} \left(\int f_{k,\ell}(|b|,t)a^{it}dt\right)v^kw^\ell$$
where each $f_{k,\ell}$ is a continuous complex function with compact 
support in $\R^+\times \R$ and such that $f_{k,\ell}(0,t)=0$ for all 
$t$, except when $k=\ell=0$. Then $A_0$ is a $^*$-algebra, acting 
non-degenerately on $\Cal H$.
\einspr

We will, as before, 
denote by $A$ the norm closure of $A_0$ and we will use $M$ for 
the weak closure. Again, both the C$^*$-algebra $A$ and the von Neumann 
algebra do not depend on the particular representation of this triple.
\snl
A simple calculation shows that, when $x$ is of the form above, then 
$$x^*=\sum_{k,\ell=0}^{1} \left(\int \overline{f_{k,\ell}}(e^{-t\theta}|b|,t)
a^{-it}dt\right) (-1)^{k\ell}v^{-k}w^{-\ell}$$
and
$$x^*x=\sum_{k,\ell,k',\ell'} \left(\iint 
\overline{f_{k,\ell}}(e^{-t\theta}|b|,t)
f_{k',\ell'}(e^{-t\theta}|b|,s)a^{i(s-t)}\,ds\,dt\right)(-1)^{\ell(k-k')}
v^{k'-k}w^{\ell'-\ell}.$$
\snl
The comultiplication $\Phi$ on this C$^*$-algebra $A$ and on this von Neumann 
algebra $M$ is implemented by a multiplicative unitary $W$ of the form
$$F(\hat b\ot b,\hat w\ot w)\exp \frac {i}{\theta}(\log \hat a\ot \log 
a)$$
where $F$ is the modified quantum exponential function (see [W-Z]). The 
choice of the dual triple is as follows:
$$\hat a=|b|^{-1}\qquad \hat b= e^{\frac12 i\theta} b^{-1}a \qquad \hat 
w=\alpha w$$
where $\alpha$ is $\pm 1$. There is some very peculiar fact however. 
Whereas the operator $W$ can be defined as a unitary for any other 
admissible triple $(\hat a,\hat b,\hat w)$, the Pentagon equation will 
not even be valid in all cases when this triple is chosen as above. 
There is a restriction on $\theta$. It must be of the form 
$\frac{\pi}{2k+3}$ with $k=0,1,2,\dots$. And then, $\alpha$ must be 
taken to be $(-1)^k$. For details, see section 2 of [W-Z].
\snl
We will not need to know more about these facts here. Just as in the 
other cases however, we need the formulas for the scaling group. They 
are more or less obvious in the case of $a$ and $b$, but not for the 
action on $w$. Again, we must refer to [W-Z].

\iinspr{5.10} Proposition \rm 
There exists a strongly continuous one-parameter group of 
$^*$-automor\-phisms $(\tau_t)$ of $A$ defined by $\tau_t(a)=a$, 
$\tau_t(b)=e^{-t\theta}b$ and $\tau_t(w)=w$. There is also an 
involutive $^*$-anti-automorphism $R$ of $A$ given by $R(a)=a^{-1}$, 
$R(b)=-e^{-\frac12 i\theta}a^{-1}b$ and $R(w)=-\alpha w$. Together, they 
give the polar decomposition of the antipode (as in 3.12).
\einspr

The first part of this proposition is not hard to prove. In fact, 
$\tau_t$ is implemented by $a^{-it}$. The action of $R$ on $a$ and 
$|b|$ is similar to the corresponding part in the complex $az+b$-case 
(cf.\ lemma 3.14).
\snl
Now we come to the right Haar measure. We have the following result.

\iinspr{5.11} Proposition \rm
There exists a lower semi-continuous densely defined faithful 
KMS-weight $\psi$ on $A$ given by
$$\psi(x)=\int f_{0,0}(u,0) du$$
when $x$ is of the form
$$x=\sum  v^kw^\ell\int f_{k,\ell}(|b|,t)a^{it}dt $$
as in proposition 5.9.
\einspr

Observe the similarity of this formula with the one in proposition 4.1.
\snl
We have 
$$\align
\psi(x^*x)&=\sum_{k,\ell} \iint |f_{k,\ell}(e^{-t\theta}u,t)|^2 du\,dt \\
          &=\sum_{k,\ell} \iint |f_{k,\ell}(u,t)|^2 e^{t\theta}du\,dt
\endalign$$
when $x$ is as above. Therefore, the GNS-space is 
$L^2(\R^+) \ot L^2(\R) \ot \C^2 \ot \C^2$ and the canonical map 
$\eta_\psi$ is given by 
$$\eta_\psi(x)=\xi=\sum \xi_{k,\ell} \ot e_k \ot e_\ell$$
where
$$\xi_{k,\ell}(u,t)=e^{\frac12 t\theta} f_{k,\ell}(u,t).$$
Again, the calculation of the different operators gives 
$$\align
\pi(a)&=a_0 \ot a_1 \ot 1 \ot 1 \\  
\pi(b)&=b_0 \ot 1   \ot s \ot 1 \\  
\pi(w)&=1   \ot 1   \ot m \ot s 
\endalign$$
where $a_0, b_0, m$ and $s$ are as in proposition 5.8 and $a_1$ is 
defined on $L^2(\R)$ by $(a_1^{is}g)(t)=g(t-s)$. For the operators 
coming from right multiplication, we get
$$\align
\pi'(a)&=e^{-\frac12 i\theta}(1 \ot a_1 \ot 1 \ot 1) \\  
\pi'(b)&=b_0 \ot c_1   \ot s \ot m \\  
\pi'(w)&=1   \ot 1   \ot 1 \ot s 
\endalign$$
where $c_1$ is defined on $L^2(\R)$ by $(c_1g)(t)=e^{t\theta}g(t)$. 
Observe that $\pi'(a)$ is not self-adjoint here.
\snl
The involution $T$ is given by $J|T|$ where 
$$\align
|T|&=1 \ot c_1^{-\frac12} \ot 1 \ot 1 \\  
J  &=\quad J_{01} \quad \ot \quad J_{23}
\endalign$$
with $J_{01}$ acting on $L^2(\R^+\times \R)$ and $J_{23}$ on $\C^2 \ot 
\C^2$ as 
$$\align 
(J_{01}\xi)(u,t)&=e^{\frac12 t\theta}\overline\xi(e^{t\theta}u,-t) \\
J_{23}(e_k\ot e_\ell)&=(-1)^ke_k\ot e_\ell.
\endalign$$  
Observe again that $|T|$ commutes with $\pi(b)$ and $\pi(w)$ but not 
with $\pi(a)$. In fact we have that 
$|T|^{it}\pi(a)|T|^{-2it}=\sigma(\pi(a))=e^{-t\theta}\pi(a)$. 
This explains why $\pi'(a)$ is not self-adjoint.
\snl
The following result will, as before, imply the right invariance of 
$\psi$. 

\iinspr{5.12} Proposition \rm Let $(a,b,w)$ be an admissible triple as in 5.8 
and let $(\hat a,\hat b,\hat w)$ be the associated triple given as 
above by $\hat a=|b|^{-1}$, $\hat b=e^{\frac12 i\theta}b^{-1}a$ and 
$\hat w=\alpha w$. Then $\hat \tau$ coincides with $\sigma$.

\snl\bf Proof: \rm 
$$\align 
\sigma_t(\hat a) & = \sigma_t(|b|^{-1})=|b|^{-1}=\hat a 
   =\hat\tau_t(\hat a)\\
\sigma_t(\hat b) & = e^{\frac12 i\theta}\sigma_t(b^{-1}a)=
   e^{\frac12 i\theta}e^{-t\theta}b^{-1}a=e^{-t\theta}\hat 
b=\hat\tau_t(\hat b)\\
\sigma_t(\hat w) & = \sigma_t(\alpha w)=\alpha w=\hat w=\hat\tau_t(\hat w).
\endalign$$
\einspr

In this case, the Haar weight is again not invariant, but only relatively 
invariant with respect to the scaling group. 
Indeed, with the notations as before, we have
$$\align
\psi(\tau_s(x))&=\psi\left(\sum_{k,\ell} v^kw^\ell
\int f_{k,\ell}(e^{-s\theta}|b|,t)a^{it}dt \right)\\
&=\int f_{0,0}(e^{-s\theta}u,0)du =e^{s\theta}\int f_{0,0}(u,0)du\\
&=e^{s\theta}\psi(x).
\endalign$$
As before, also here the regular representation $\widetilde W$ can be 
calculated using the exponential properties of the (quantum) 
exponential function involved. We find
$$\widetilde W=F(\hat b\ot b)\exp\frac{i}{\theta}(\log \hat a\ot \log 
a)$$
where now
$$\align 
a &= a_0 \ot a_1 \ot 1 \ot  1\\
b &= b_0 \ot 1   \ot s \ot  1\\
w &= 1   \ot 1   \ot m \ot  s
\endalign$$
and
$$\align 
\hat a &= 1 \ot c_1 \ot 1 \ot 1\\
\hat b &= \cdot \ot a_1 \ot s \ot \cdot \\
\hat w &= -\alpha (1 \ot 1 \ot m \ot 1)
\endalign$$
and where the part $\cdot \ot \cdot$ as acting on $L^2(\R^+)\ot \C^2$ is 
given by 
$$e^{\frac12 i\theta}b_0^{-1}a_0 \ot 1 - b_0^{-1} \ot m$$
\nl	 
The {\it last example} we consider is the quantum $E(2)$ as introduced by 
Woronowicz. 
\snl
The situation with this example is, in several ways, different from the 
previous ones. The main difference is that the quantum $E(2)$ is not 
self-dual. This implies that the strategy to determine the Haar measure 
must be modified. However, this is quite interesting as we come closer 
to the general setting which we will discuss briefly in the next section. The 
other difference is that this example has been already studied some 
time ago by Woronowicz. In particular, the multiplicative unitary is 
only touched and manageability was not yet known. On the other hand, 
the Haar measures have been obtained already by S.\ Baaj in [B1] and 
[B2], but 
the treatment seems to be more complicated than ours. For all these 
reasons, here we will only briefly indicate how our method is applied in 
this case and we will give more details in a separate paper that we 
plan to write [J-VD]. There, we will show that the quantum $E(2)$ is 
indeed a locally compact quantum group in the sense of [K-V2]. We will 
also use the opportunity to revise and update the treatment of this 
example.  
\snl
As before, we first describe the operators involved.

\iinspr{5.13}  Definition \rm
Let $q\in \R$ and suppose $0<q<1$ as before. Consider $\ell^2(\Z)$ with 
an orthonormal basis $\{e_k \mid k\in \Z \}$. Define a unitary operator 
$s$  and a non-singular positive self-adjoint operator $m$ on $\ell^2(\Z)$
by 
$$se_k=e_{k+1} \qquad \text{ and } \qquad me_k=q^ke_k.$$
Then consider $\Cal H=\ell^2(\Z) \otimes \ell^2(\Z)$ and define 
operators $a, b, c$ and $d$ on $\Cal H$ by
$$\alignat 2 a&=m^{-\frac12}\otimes m &\qquad\qquad  b&=m^{\frac12} \otimes s \\
          c&=s \otimes s &\qquad\qquad  d&=s\otimes m^{-1}.
\endalignat$$
\einspr

The operators $s$ and $m$ are the same as in 5.2 and so $ms=qsm$ as 
before. The operators $b$ and $d$ are normal operators. The operator 
$c$ is unitary and the operator $a$ is again non-singular, positive 
self-adjoint. We have the commutation rules $cd=qdc$ and $ab=qba$ 
(where as usual, the last one is interpreted as 
$a^{it}ba^{-it}=q^{it}b$ for all $t\in \R$). This means that the pair 
$(c,d)$ 'generates' the quantum group $E_q(2)$ while the pair $(a,b)$ 
generates its dual $\hat E_q(2)$ as in [VD-W]. Observe that in [VD-W], $c$ 
is denoted by $v$ and $d$ is denoted by $n$. Moreover, $\mu$ is used 
instead of $q$. Finally, $a=\mu^{\frac12} N$. 
\snl
With these notations, it is easy to verify the following.

\iinspr{5.14} Proposition  \rm
The operators $a,b,c$ and $d$ satisfy the commutation rules of section 5 of 
[W4].
\snl\bf Proof: \rm
We will only consider the non-trivial condition vi) of section 5 in [W4]. 
With our 
notations, we need to have
$$db-q^{\frac12}bd=(1-q^2)q^{-\frac12}a^{-1}c.$$
Now
$$\align db & = s m^{\frac12} \otimes m^{-1} s= q^{-\frac12}
(m^{\frac12} s \otimes m^{-1} s)= q^{-\frac12} a^{-1}c\\
db & = m^{\frac12} s \otimes s m^{-1} = q (m^{\frac12} s \otimes m^{-1} 
s) = q a^{-1} c
\endalign
$$
and the equation above follows easily.
\einspr

Of course, we should be more careful with the above relations of 
unbounded operators. Observe that also Woronowicz is only formal when 
writing down his formula in [W4]. In [J-VD], we plan to be more precise.
\snl
Then, as Woronowicz claims, we get the following result.

\iinspr{5.15} Proposition \rm
Define $W=F(ab \otimes cd)\chi(a\otimes 1, 1\otimes c)$ where 
$\chi(q^{\frac12 n},z)=z^n$ for $n\in \Z$ and $z\in \C$ with $|z|=1$ 
and where $F$ is a function defined on the appropriate part of $\C$ by
$$F(z)=\prod_{k=0}^\infty \frac{1+q^{2k}\bar z}{1+q^{2k}z}. $$
Then $W$ satisfies the Pentagon equation.
\einspr

Observe that $\chi(a\otimes 1, 1\otimes c)$ is nothing else but 
$(1\otimes c)^{(N\otimes 1)}$ in the notation of Woronowicz. Moreover, 
$ab\otimes cd$ is a normal operator when defined in the appropriate 
way. The function $F$ is first only defined whenever the denominator in 
the infinite product is non-zero. Then, it is extended continuously 
along concentric circles. The values of $F$ are again in the 
unit circle and so the operator $W$ is unitary.
\snl
We know that formally, the comultiplication $\hat \Delta$ on $\hat 
E_q(2)$ is given by the formulas
$$\align \hat \Delta(a) & = a\otimes a \\
\hat \Delta(b) & = a \otimes b + b \otimes a^{-1}
\endalign
$$
(see e.g.\ [W4]). 
The antipode 
$\hat S$ is given by $\hat S(a)=a^{-1}$ and $\hat S(b)=-q^{-1}b$. We 
will show in [J-VD] that, as expected from the general theory, $(\hat 
S\otimes \iota)W = W^*$.
\snl
The right leg of $W$ gives the C$^*$-algebra $A$ (and the von Neumann 
algebra $M$), 'generated' by the elements $c$ and $d$ (cf.\ e.g.\ [W4]). 
And of course, $W$ induces the comultiplication on $A$ (and on $M$) in 
the usual way.
\snl
Now, it is relatively easy to prove the main result concerning this 
example.

\iinspr{5.16} Theorem \rm 
The right invariant Haar weight on the quantum group $E_q(2)$ is given 
by the formula
$$\psi(x)=\sum_{k=-\infty}^{+\infty} q^{-2k}\langle x(e_0\otimes 
e_k),e_0\otimes e_k\rangle$$
whenever $x\in A^+$.
\einspr

We will not give details of the proof here but refer to [J-VD]. However, 
the idea behind it is as follows and is completely similar to the 
previous situations.
\snl
First consider 
$$\psi_0(x)=\sum_{k,\ell=-\infty}^{+\infty} q^{-2k}\langle x(e_\ell\otimes 
e_k),e_\ell\otimes e_k\rangle$$
for any positive operator $x$ on $\Cal H$. We have $\psi_0=\text{Tr} 
(h\,\cdot)$ where $h=1\otimes m^{-2}$. This is a faithful normal 
semi-finite weight on $\Cal B(\Cal H)$. The modular automorphisms 
$\sigma_t$ are implemented by $h^{it}=1\otimes m^{-2it}$. If we 
'restrict' this to the operators $a$ and $b$, we get $\sigma_t(a)=a$ 
and $\sigma_t(b)=q^{-2it} b$ for all $t$. On the other hand we had 
$\hat S^2(a)=a$ and $\hat S^2(b)= q^{-2}b$. So we see that, also here, 
the modular automorphism group $\sigma$ coincides with the scaling 
group $\hat \tau$. This is precisely what is needed for right 
invariance. However, the weight $\psi_0$ is not semi-finite on $M$. 
Fortunately, we can cut it down to $\psi=\psi_0((p\otimes 1) \,\cdot\, 
(p\otimes 1))$ where $p$ is the one-dimensional projection  on the space 
$\C e_0$. This will give a semi-finite weight on our von Neumann 
algebra $M$. The fact that $p\otimes 1$ commutes with both $a$ and $b$ 
guarantees that the argument to prove invariance can still be used in 
this situation.
\snl
In [J-VD] we not only plan to give the details of the proof, but we will 
also describe the right regular representation and compare it with the 
formula obtained by S.\ Baaj in [B1] and [B2].

\nl\nl	  


\bf 6. Conclusions and perspectives \rm
\nl
In this paper, we constructed the Haar measures for certain locally 
compact quantum group candidates. We have done this in fairly great 
detail for the quantization of the $az+b$-group with a complex 
deformation parameter. We have also treated other cases: the 
quantization of the $az+b$-group with a real parameter, the quantum 
$ax+b$-group and (very briefly) the quantum $E(2)$. We also constructed 
regular representations.
\snl
Now we would like to explain why the method that we have used should 
always lead to a solution. To see this, start with a locally compact 
quantum group in the sense of Kustermans and Vaes [K-V2] and consider 
the von Neumann algebra setting [K-V4]. Denote by $M$ the underlying 
von Neumann algebra and by $\hat M$ the dual von Neumann algebra, both 
acting on the GNS-space $\Cal H$ of the Haar weights. Use $\psi$ to 
denote the right invariant Haar weight on $M$.
\snl
In [VD6] we will show that there is a faithful semi-finite normal 
weight $f$ on $\Cal B(\Cal H)$ such that $x\mapsto f(y^*xy)$ is a 
scalar multiple of $\psi$ for certain well-chosen elements $y$ in the 
commutant of $\hat M$. This weight is of the form 
$\text{Tr}(h\,\cdot\,)$ where $h$ is a certain implementation of the 
analytic generator of the scaling group. In fact, it is the positive 
operator that realizes the manageability of the multiplicative unitary. 
The modular automorphism group of $f$ will coincide with the scaling 
group on $\hat M$. On the other hand, if $f$ is a weight with this last 
property, it has to be the above one. Therefore, $x\mapsto f(y^*xy)$ 
will be right invariant. And there must be elements $y$ for which this 
gives a semi-finite weight, hence the right Haar measure.
\snl
This property is behind our constructions in this paper. One must be 
aware of the fact that the multiplicative unitaries we start with in 
all cases are not the regular representations (as we have mentioned 
already). In all the examples, except for the quantum $E(2)$ in section 
5, the left and the right leg generated the same von Neumann algebra, 
namely the full $\Cal B(\Cal H)$. So, the commutant is trivial and we 
did not need to cut down as we would have to do when the multiplicative 
unitary was the regular representation. In the case of $E(2)$, we are 
in an intermediate situation and we do have to cut down using an 
element in the commutant (see the remark after theorem 5.16).
\snl
All of this will be explained in detail in [VD6]. In [VD5] we will do 
this in a purely algebraic context (for multiplier Hopf algebras with 
integrals, the so-called algebraic quantum groups). This paper will 
give the full algebraic background and will help to understand what 
will be done in [VD5].
\snl
A theory of locally compact quantum groups with a set of natural 
axioms, not assuming but proving the existence of the Haar measures, 
still seems to be out of reach. On the other hand, what we claim here 
is that, whenever the Haar measure exists, our method should provide 
it. This might be even more important when constructing examples than a 
theoretical existence proof. And the results obtained here and techniques 
that are used could contribute to such a theory where the existence of 
the Haar measures is proved from the axioms.
  
\nl \nl


\bf Appendix: Heisenberg commutation relations.  \rm
\nl
At various places in this paper, we encounter a pairing between abelian 
groups, identifying one as the dual of the other, and a compatible set 
of representations of these two groups. In this appendix, we will first 
discuss the uniqueness property of such a pair of representations in 
the general case. Then we will be more concrete and look at the various
cases.
\snl
The starting point is a pair of abelian locally compact qroups $G$ and 
$K$ and a non-degenerate continous pairing $\langle\ ,\ \rangle$ from 
$G\times K$ to the unit circle $\T$. We assume that this pairing is 
multiplicative in both variables so that in fact, we have a bicharacter. 
Each of the two groups can be identified with the Pontryagin dual of 
the other one.

\inspr{A.1} Definition \rm 
A pair of continuous unitary representations $\pi$ of $G$ and $\gamma$ 
of $K$ satisfy the {\it Heisenberg commutation relations} if
$$\gamma(k)\pi(g)=\langle g,k\rangle \pi(g)\gamma(k)$$
for all $g\in G$ and $k\in K$. 
\einspr

The typical example is the following.

\inspr{A.2} Example \rm
Let $G$ be a locally compact abelian group, let $K$ be the 
Pontryagin dual $\hat G$ and write $\langle g,p\rangle$ for the value of the 
element $p\in \hat G$ in the point $g\in G$. This is indeed a 
non-degenerate bicharacter on $G\times \hat G$. Consider the Hilbert 
space $L^2(G)$ where $G$ is considered with its Haar measure. Define 
representations $\pi_0$ of $G$ and $\gamma_0$ of $\hat G$ by
$$\align (\pi_0(g)\xi)(h)&=\xi(g^{-1}h)\\
        (\gamma_0(p)\xi)(h)&=\langle h, p\rangle\xi(h)
\endalign$$	
where $\xi \in L^2(G)$, $g,h\in G$ and $p\in \hat G$. A simple 
calculation gives the Heisenberg commutation rules as in the definition 
above.
\einspr

The main result about such a Heisenberg representation is the 
following.

\inspr{A.3} Proposition \rm
Given a pair of two locally compact abelian groups as above, then there 
is (up to unitary isomorphism) just one irreducible Heisenberg 
representation. Moreover, any Heisenberg representation is unitary 
equivalent with a multiple of this irreducible representation.
\einspr

More precisely (or equivalently), given a locally compact abelian group 
$G$ with dual group $\hat G$ and continuous unitary representations 
$\pi$ of $G$ and $\gamma$ of $\hat G$ on a Hilbert space $\Cal H$ 
satisfying
$$\gamma(p)\pi(g)=\langle g,p\rangle \pi(g)\gamma(p)$$
for all $g\in G$ and $p\in \hat G$, there exists a Hilbert space $\Cal 
K$ and a unitary $U:L^2(G)\otimes \Cal K \to \Cal H$ such that
$$\align \pi(g)&=U(\pi_0(g)\otimes 1)U^* \\
         \gamma(p)&=U(\gamma_0(p)\otimes 1)U^*
\endalign$$
for all $g\in G$ and $p\in \hat G$. Here, $\pi_0$ and $\gamma_0$ are 
the representations given in the example A2 above.
\snl
The result is certainly well known (see e.g.\ [??]). Nevertheless, for 
completeness, let us give some possible argument.
\snl
The representation $\gamma$ of $\hat G$ gives a non-degenerate 
representation of the C$^*$-algebra $C_0(G)$ of complex continuous 
functions on $G$ tending to $0$ at infinity. The group $G$ acts on 
$C_0(G)$ by left translation. Denote this action by $\alpha$. The 
Heisenberg commutation rules guarantee that we have a covariant 
representation of this covariant system $(C_0(G),G,\alpha)$. This gives 
rise to a non-degenerate representation of the crossed product 
$C_0(G)\times_\alpha G$. In fact, any non-degenerate representation of 
this crossed product will come from a Heisenberg representation of the 
pair $(G,\hat G)$ in this way.
\snl
Because the group is abelian, the full crossed product coincides with 
the reduced crossed product. This last one can be realized on the 
tensor product $L^2(G)\otimes L^2(G)$ coming from the 
Heisenberg pair $(\pi_1,\gamma_1)$ given by
$$\align 
\pi_1(g)&=\pi_0(g)\otimes \pi_0(g)\\
\gamma_1(p)&=\gamma_0(p)\otimes 1
\endalign$$
where $(\pi_0,\gamma_0)$ is as before and $g\in G$ and $p\in \hat G$.
The unitary $W$ on $L^2(G)\otimes L^2(G)$ defined by 
$(W\xi)(g,h)=\xi(g,g^{-1}h)$ gives
$$\align 
\pi_1(g)&=W(\pi_0(g)\otimes 1)W^* \\
\gamma_1(p)& =W(\gamma_0(p)\otimes 1)W^*.
\endalign$$
This can be used to show that the crossed product $C_0(G)\times_\alpha 
G$ is isomorphic with the C$^*$-algebra of the compact operators on 
$L^2(G)$. Then, it follows from a general result about the 
representation theory of compact operators (see e.g.\ [??]) that any 
non-degenerate representation of this crossed product is equivalent 
with a multiple of the unique irreducible representation. All of this will 
imply the result of proposition A.3.
\nl
We now collect the different special cases that we consider in this 
paper.

\inspr{A.4} Examples \rm
i) First we have the finite cyclic group $\Z_{2n}$. We use the pairing 
with itself given by $\langle k,\ell \rangle=q^{k\ell}$ where 
$q=\exp{\frac{\pi i}{n}}$. The irreducible Heisenberg representation is 
given in proposition 3.2 when we take $s$ to be the generator of the 
first factor $\Z_{2n}$ and $m$ as the generator of the second factor 
$\Z_{2n}$ in the pairing.
\snl
ii) The second case is the group $\R$ paired with itself by using 
$\langle t,s\rangle=\exp \frac{\pi i ts}{n}$. The irreducible 
Heisenberg representation given in proposition 3.4 is not the standard 
one of example A.2.  We have a representation acting on $L^2(\R^+)$ 
defined by 
$$\pi(t)=a_0^{it}\qquad\text{ and }\qquad \gamma(s)=b_0^{is}.$$
\snl
iii) In definition 3.1 we are essentially dealing with the self-dual group 
$\Z_{2n}\times \R$ (which is isomorphic with the group $\Gamma$ as 
introduced in 3.2 ii).
\snl
iv) In section 4 we have another Heisenberg representation of the 
pairing in ii) above. It is essentially the standard one, acting on 
$L^2(\Bbb R)$, and given by   
$$\pi(t)=c_1^{it}\qquad\text{ and }\qquad \gamma(s)=a_1^{is}.$$
\snl
v) In the first example of section 5, we have the self-dual group 
$\T\times \Z$. The pairing is coming from the pairing between $\Z$ and 
$\T$ which is given by $\langle n,z\rangle= z^n$ for $n\in \Z$ and $z\in 
\T$. The irreducible representation we used is 
$$\pi(n)=s^n \qquad\text{ and }\qquad \gamma(q^{it})=m^{it}$$
with $n\in \Z$ and $t\in \R$.
The irreducible representation for the product $\T\times \Z$, paired 
with $\Z\times \T$, is essentially the tensor product of two copies 
above.
\snl
vi) In the second example of section 5 we have the group $\Z_2\times 
\R$ with the usual irreducible representation for $\Z_2$ (as for 
$\Z_{2n}$ in example i)) and the irreducible representation for $\R$ 
more or less as in ii) and in iii).
\snl
vii) Finally, in the last example of section 5, the building block is again the same 
as in example v) for the natural pairing between $\Z$ and $\T$.
\einspr

We finish this appendix with formulating a result, related with the 
material here  and also used in this paper (cf.\ lemma 3.14). 

\inspr{A.5} Proposition \rm
Let $(a,b)$ be a pair of positive, non-singular self-adjoint 
operators on a Hilbert space such that $a^{it}ba^{-it}=e^{-\frac{\pi 
t}{n}}b$ (where $n\in \N$). Then $(e^{-\frac{\pi i 
t^2}{2n}}a^{-it}b^{it})_{t\in\R}$ is again a strongly continuous 
one-parameter group of unitaries whose analytic generator is the 
closure of $e^{\frac{\pi i}{2n}}a^{-1}b$.
\einspr

The result is also standard (see e.g.\ [??]). To verify that these 
unitaries form a one-parameter group, just observe
$$\align e^{-\frac{\pi i (t+s)^2}{2n}} a^{-it} a^{-is} b^{it} b^{is} 
&=e^{-\frac{\pi i (t+s)^2}{2n}}
e^{\frac{\pi i ts}{n}} a^{-it}b^{it}a^{-is}b^{is} \\
&=(e^{-\frac{\pi i t^2}{2n}} a^{-it} b^{it})
(e^{-\frac{\pi i s^2}{2n}} a^{-is}b^{is}).
\endalign$$
To prove that the generator is the closure of $e^{\frac{\pi 
i}{n}}a^{-1}b$ is of course more difficult. The following formal 
calculation should at least clarify why this is true:
$$\align (e^{\frac{\pi i}{2n}} a^{-1}b)^p
&= e^{\frac{\pi i p}{2n}} (a^{-1}b a^{-1}b \dots a^{-1}b)\\
&= e^{\frac{\pi i p}{2n}} e^{\frac12 p(p-1)\frac{\pi i}{n}\,} a^{-p}b^p\\ 
&= e^{\frac{\pi i p^2}{2n}} a^{-p}b^p
\endalign$$
where we used the basic commutation rule $ba^{-1}=e^{\frac{\pi 
i}{n}}a^{-1}b$.
\snl
Observe that $(a,b)$ satisfies the same commutation rule as 
$(a,e^{\frac{\pi i}{2n}}a^{-1}b)$. By the above theory we know that 
(up to a possible multiplicity), these pairs are unitarily equivalent. 
In fact, they are unitarily equivalent. The unitary relating them was 
given in lemma 3.14.
\nl\nl


\parindent=20pt
\noindent
\bf References \rm
\bigskip

\item{[A]} E.\ Abe: \it Hopf algebras. \rm Cambridge University Press (1977).
\smallskip

\item{[Ar]} W.\ Arveson: \it An invitation to C$^*$-algebras. \rm
Springer-Verlag, New York (1976).
\smallskip

\item{[B1]} S. Baaj: \it Repr\'{e}sentation r\'{e}guli\`{e}re du groupe 
quantique $E_\mu(2)$ de Woronowicz.
 \rm 
C.R. Acad. Sci. Paris, S\'er. I {\bf 314} 
(1992) 1021-1026.
\smallskip

\item{[B2]} S. Baaj: \it Repr\'{e}sentation r\'{e}guli\`{e}re du groupe 
quantique des d\'{e}placements de Woronowicz.
 \rm
Ast\'erisque {\bf 232} (1995) 11-48.
\smallskip

\item{[B-S]} S.\ Baaj \& G.\ Skandalis: \it Unitaires multiplicatifs et dualit\'e
          pour les produits crois\'es de C$^*$-alg\`ebres. \rm Ann.
          Scient.\ Ec.\ Norm.\ Sup., 4\`eme s\'erie, {\bf 26} (1993) 425-488.
\smallskip

\item{[D]} V.G.\ Drinfel'd: \it Quantum groups. \rm Proceedings ICM Berkeley
          (1986) 798-820.
\smallskip

\item{[E-S]} M.\ Enock \& J.-M.\ Schwartz: {\it Kac algebras and duality
for locally compact groups.} Springer (1992).
\smallskip

\item{[J-VD]}  A.\ Jacobs \& A.\ Van Daele: {\it The quantum $E(2)$ as 
a locally compact quantum group.} Preprint K.U.\ Leuven (in 
preparation).
\smallskip
 
\item{[J]} Jimbo: {\it A q-analogue of $U(g)$ and the Yang-Baxter 
equation}. Lett.\ Math.\ Phys.\ {\bf 10} (1985) 63-69.
\smallskip

\item{[K-R]} R.V.\ Kadison \& R.\ Ringrose: \it Fundamentals of the theory of
operator algebras. \rm Academic Press, Orlando (1986).

\smallskip

\item{[Ki]} E.\ Kirchberg: Lecture at the conference 'Invariants in 
operator algebras'. Copenhagen (1992).
\smallskip

\item{[K-S]} A.\ Klimyk \& K.\ Schm\"udgen: {\it Quantum groups and 
their representation}. Springer (New York) 1997.
  
\item{[K1]} J.\ Kustermans: {\it C$^*$-algebraic quantum groups arising
from algebraic quantum groups.} Ph.D.\ thesis K.U.\ Leuven (1997).
\smallskip

\item{[K2]} J.\ Kustermans: {\it The analytic structure of algebraic 
quantum groups.} Preprint K.U.\ Leuven (2000) (Funct-An/970710).
\smallskip

\item{[K-V1]} J.\ Kustermans \& S.\ Vaes: \it A simple
definition  for locally compact quantum groups. \rm  C.R. Acad. Sci., Paris,
S\'er. I {\bf 328 (10)} (1999) 871-876.

\smallskip

\item{[K-V2]} J.\ Kustermans \& S.\ Vaes: \it Locally compact
quantum groups. \rm Ann.\ Sci.\
 Ec.\ Norm.\ Sup.\ {\bf 33} 
(2000), 837-934.
\smallskip

\item{[K-V3]} J.\ Kustermans \& S.\ Vaes: \it The operator algebra
approach to quantum groups. \rm Proc. Natl. Acad. Sci. USA {\bf 97 (2)}
(2000) 547-552.

\smallskip

\item{[K-V4]} J.\ Kustermans \& S.\ Vaes: \it Locally quantum
groups in the von Neumann algebra setting.  \rm Preprint
K.U.Leuven (2000).

\smallskip

\item{[K-VD]} J.\ Kustermans \& A. Van Daele:
{\it C$^*$-algebraic quantum groups arising
from algebraic quantum groups.} Int. J. Math. {\bf 8} (1997) 1067-1139.
\smallskip

\item{[M-N]} M.\ Masuda \& Y.\ Nakagami: {\it A von Neumann algebra
framework for the duality of quantum groups.} Publ. RIMS Kyoto
{\bf 30} (1994) 799-850.
\smallskip

\item{[M-N-W]} M.\ Masuda, Y.\ Nakagami \& S.\ Woronowicz (in
preparation).
 Lectures at the Fields Institute and at the University of 
Warsaw (1995).
\smallskip

\item{[P]} G.K.\ Pedersen: \it C$^*$-algebras and their automorphism
groups. \rm Academic Press (1979).
\smallskip

\item{[Sa]} S.\ Sakai: \it C$^*$-algebras and W$^*$-algebras. \rm Springer
        Verlag (1971).
\smallskip

\item{[S-Z]} S.\ Stratila \& L.\ Zsid\'{o}: {\it Lectures on von Neumann 
algebras}.  Abacus Press, Tunbridge Wells, England (1979).

\smallskip

\item{[St]} S.\ Stratila: {\it Modular theory in operator algebras.} 
Abacus Press, Tunbridge Wells, England (1981).
\smallskip

\item{[Sw]} M.E.\ Sweedler: \it Hopf algebras. \rm Mathematical Lecture Note
          Series. Benjamin (1969).
\smallskip

\item{[V]} S.\ Vaes: {\it Examples of locally compact quantum groups 
through the bicrossed product construction.} To appear in the 
Proceedings of the XIIIth International Conference Mathematical Physics, 
London 2000.
\smallskip

\item{[V-V]} S.\ Vaes \& L.\ Vaynerman: {\it Extensions of locally 
compact quantum groups and the bicrossed product construction.} 
Preprint Max Planck Institut f\"ur Mathematik MPI 2001-2 (2001). 
\smallskip

\item{[V-VD]} S.\ Vaes \& A.\ Van Daele: {\it Hopf C$^*$-algebras}. 
Proc.\ London Math.\ Soc.\ {\bf 82} (2001) 337-384.
\smallskip

\item{[VD1]} A.\ Van Daele: \it Dual pairs of Hopf $^*$-algebras. \rm
Bull.\ London Math.\ Soc.\ {\bf 25} (1993) 209-230.
\smallskip

\item{[VD2]} A.\ Van Daele: {\it Multiplier Hopf algebras.} Trans.\
    Amer.\ Math.\ Soc.\ {\bf 342} (1994) 917-932.
\smallskip


\item{[VD3]} A.\ Van Daele: {\it An algebraic framework for group duality.}
Adv. in Math. {\bf 140} (1998) 323-366.
\smallskip

\item{[VD4]} A.\ Van Daele: {\it A dual pair approach to some locally 
compact quantum groups.} In preparation.
\smallskip

\item{[VD5]} A.\ Van Daele: {\it The Heisenberg commutation relations 
for an algebraic quantum group.} (I and II) In preparation (with J.\ 
Kustermans).
\smallskip

\item{[VD6]} A.\ Van Daele: {\it The Heisenberg commutation relations, 
commuting squares and the Haar measure on locally compact quantum 
groups.} In preparation (with S.\ Vaes).
\smallskip

\item{[VD-W]} A.\ Van Daele \& S.L.\ Woronowicz: {\it Duality for the 
quantum $E(2)$ group.} Pac. J. Math. {\bf 7} (1996) 255-263. 
\smallskip

\item{[W1]} S.L.\ Woronowicz: {\it Compact Matrix Pseudogroups.}
Commun. Math. Phys. {\bf 111} (1987) 613-665.

\smallskip

\item{[W2]} S.L.\ Woronowicz: {\it  Compact Quantum Groups.}
Quantum symmetries/Symm\'{e}tries quantiques.  Proceedings of the
Les Houches summer school 1995, North-Holland, Amsterdam (1998),
845--884.
\smallskip

\item{[W3]} S.L.\ Woronowicz: {\it From multiplicative unitaries to quantum 
groups.} Int. J. Math. {\bf 7} (1996) 127-149.

\smallskip

\item{[W4]} S.L.\ Woronowicz: {\it Quantum $E(2)$ and its Pontryagin 
dual.} Lett. Math. Phys. {\bf 23} (1991) 251-263.
\smallskip

\item{[W5]} S.L.\ Woronowicz: {\it Quantum $az+b$-group on complex 
plane.} Preprint University of Warsaw and University of Trondheim 
(2000). To appear in International Journal of Mathematics.
\smallskip

\item{[W6]} S.L.\ Woronowicz: {\it Quantum exponential function}. 
Reviews in Mathematical Physics. {\bf 12} (2000) 873-920. 
\smallskip

\item{[W-Z]} S.L.\ Woronowicz \& S.\ Zakrzewski: 
{\it Quantum $ax+b$-group}. Preprint University of Warsaw (2001). 
\smallskip

\end